\documentclass{amsart}
\usepackage{amsmath,amsfonts,latexsym,amssymb}

\newtheorem{theorem}{Theorem}
\newtheorem{lemma}{Lemma}

\newtheorem{proposition}{Proposition}

\theoremstyle{remark}

\newtheorem*{ack}{Acknowledgements}

\begin{document}

\markboth{Ritabrata Munshi}{Pairs of quadrics in 11 variables}
\title[Pairs of quadrics in 11 variables]{Pairs of quadrics in 11 variables}

\author{Ritabrata Munshi}   
\address{School of Mathematics, Tata Institute of Fundamental Research, 1 Dr. Homi Bhabha Road, Colaba, Mumbai 400005, India.}     
\email{rmunshi@math.tifr.res.in}

\begin{abstract}
For non-singular intersections of pairs of quadrics in $11$ or more variables, we prove an asymptotic for the number of rational points in an expanding box.
\end{abstract}


\maketitle


\section{Introduction}
\label{intro}

In 1962 Birch \cite{B}, following Davenport \cite{D}, used the Hardy-Littlewood circle method to prove a very general statement about the asymptotic count for the number of rational points on a variety defined by a system of forms of a given degree when the dimension is suitably large. For hypersurfaces (i.e. varieties defined by one equation) the result has been improved in some cases. Notably for cubic hypersurfaces, Heath-Brown \cite{HB2} has reduced the required number of variables from 17 to 14 under certain conditions. 
But for varieties given by more than one equation, Birch's result has not been improved without imposing serious restrictions on the nature of the forms.
\\

In this paper we are concerned with non-singular varieties $V$ which are given by the intersections of two quadrics, say $Q_i(\mathbf{x})= ~^t\mathbf{x}M_i\mathbf{x}=0$, $i=1,2$, in $n$ variables. It is a non-trivial consequence of the non-singularity of $V$ (see Chapter~2 of \cite{R}) that one of the matrices $M_i$ (say $M_2$) has to be of full rank. 
The singular locus, in the sense of Birch, is given by 
$$
V^\star=\left\{\mathbf{x}\in \mathbb{C}^n: \mathrm{rank}(M_1\mathbf{x},M_2\mathbf{x})\leq 1\right\},
$$
which coincides with the union of the eigenspaces of $M_2^{-1}M_1$. Let $\Delta$ be the dimension of the singular locus, which is same as maximum dimension of the eigenspaces. Since we are assuming that $V$ is non-singular we have $\Delta=1$ or $2$.\\ 

It follows from Thorem~1 of Birch \cite{B} that for any box $\mathcal{P}$ one has
\begin{align}
\label{birch}
\mathop{\sum\dots\sum}_{\substack{\mathbf{m}\in\mathbb{Z}^n\cap B\mathcal{P}\\Q_1(\mathbf{m})=Q_2(\mathbf{m})=0}}1 =\mathfrak{S}J_0(\mathcal{P}) B^{n-4}+O\left(B^{n-4-\delta}\right)
\end{align}
where $\delta>0$ if the number of variables $n\geq 13+\Delta$. Here $\mathfrak{S}$ is the standard singular series and $J_0(\mathcal{P})$ is the singular integral. One seeks to improve this and reduce the required number of variables (at least) to $n\geq 9$. Indeed from the work of Demyanov \cite{De} it is known that for $n\geq 9$ local solutions exist for any finite prime (also see \cite{BLM}). On the other hand from the work of Colliot-Th\'el\`ene, Sansuc and Swinnerton-Dyer \cite{CSS}, \cite{CSS2}, we know that the Hasse principle holds for $n\geq 9$. (This has been improved for non-singular intersections in a recent pre-print of Heath-Brown \cite{HB4}, where he shows that the Hasse principle holds for $n\geq 8$.)\\ 

The purpose of this paper is to use a new form of the circle method (`nested $\delta$-method') to extend the admissible range for $n$ in \eqref{birch} when $\Delta=1$, which is the generic situation. Our result gives an asymptotic for the counting function for $n\geq 11$. Moreover we believe that this is the first occasion where the $\delta$-method is used to deal with more than one (non-linear) equation. \\

Let $W$ be a compactly supported non-negative smooth function on $\mathbb R^n$ whose support does not intersect the singular locus, i.e. $\text{Supp}(W)\cap V^\star=\emptyset$. Suppose that the derivatives of $W$ satisfy the bound $W^{(j)}\ll_j H^j$ where $H\geq 1$. For example $W$ can be such that $W(\mathbf{y})=1$ for $\mathbf{y}\in \prod_{i=1}^n [a_i,b_i]$ and $W(\mathbf{y})=0$ for $\mathbf{y}\notin \prod_{i=1}^n [a_i-H^{-1},b_i+H^{-1}]$. The extra parameter $H$ will be used later to obtain an asymptotic for the counting function without any smooth weight. \\  

\begin{theorem}
\label{mthm}
Let $Q_i(\mathbf{x})$, with $i=1,2$, be two quadratic forms with rational coefficients in $n$ variables with matrices $M_i$. Suppose the variety given by $Q_1(\mathbf{x})=Q_2(\mathbf{x})=0$ is non-singular. Suppose $M_2$ is of full rank and the eigenvalues of $M_2^{-1}M_1$ are all distinct (i.e. $\Delta=1$). Let $W$ be a smooth function as above. Then we have
\begin{align}
\label{new-birch}
\mathop{\sum\dots\sum}_{\substack{\mathbf{m}=(m_1,\dots,m_n)\in\mathbb{Z}^n\\Q_1(\mathbf{m})=Q_2(\mathbf{m})=0}}W\left(\frac{\mathbf{m}}{B}\right)=\mathfrak{S}J_0(W)B^{n-4}+O\left(H^2B^{n-5+\varepsilon}+H^{2n}B^{3n/4-41/32+\varepsilon}\right),
\end{align}
where the singular integral $J_0(W)$ depends on $W$, and the implied constant depends on $Q_i$ and $\varepsilon$.\\
\end{theorem}

In Section~\ref{maintermsection} we will show that the singular integral is given by
\begin{align*}
J_0(W)=\mathop{\lim\lim}_{\varepsilon_1,\varepsilon_2\rightarrow 0}\frac{1}{4\varepsilon_1\varepsilon_2}\int_{\substack{|Q_1(\mathbf{y})|<\varepsilon_1\\|Q_2(\mathbf{y})|<\varepsilon_2}}W(\mathbf{y})\mathrm{d}\mathbf{y},
\end{align*}
and the singular series is given by
$$
\mathfrak{S}=\sum_{q_1=1}^\infty\;\sum_{q_2=1}^{\infty}\;\frac{C_{q_1,q_2}(\mathbf{0})}{q_1^{n-1}q_2^{n}},
$$
where $C_{q_1,q_2}$ is as defined in \eqref{char-sum}. One can show that the above series, which is also given by an Euler product, in fact coincides with  the product of the local densities. The technical hypothesis regarding the support of the weight function $\text{Supp}(W)\cap V^\star=\emptyset$, is required in our estimation of the main term in Section~\ref{maintermsection}. One should relate this with Birch's treatment of the singular integral in Section~6 of \cite{B}, where he needs to `dig out the singular points from the box' $\mathcal{P}$. More precisely, he breaks the box  into several smaller pieces, and then evaluates the integral for each smaller box in two different ways depending on whether the box intersects the singular locus or not. Perhaps a similar trick can be used in the present scenario to get rid of the technical hypothesis.\\

Suppose $W$ is such that $H=1$ and $J_0(W)>0$. Then the error term is smaller than the main term if $n-4>3n/4-41/32$ which holds if $n\geq 11$. Note that for $\Delta=1$, Birch \cite{B} required $n\geq 14$. The asymptotic formula is the first improvement over \eqref{birch}, in such generality, for $n=11, 12$ and $13$. In special situations, however, there are few instances which go far beyond Birch's result. In the case of pairs of diagonal quadratic forms an asymptotic for the count can be obtained for fewer variables. This was done by Cook \cite{C} for diagonal quadric pairs with nine or more variables ($n\geq 9$). Recently in two joint works with T.D. Browning \cite{BM1}, \cite{BM2}, we established an asymptotic for the counting functions when the pair has a special structure, namely $Q_1(\mathbf{x})=q_1(x_1,\dots,x_{n-2})-x_{n-1}^2-x_n^2$ and $Q_2(\mathbf{x})=q_2(x_1,\dots,x_{n-2})$, with $q_1$ and $q_2$ being quadratic forms (not necessarily diagonal). Such pairs of quadrics appear naturally in many other important counting problems, e.g. Batyrev-Manin conjecture for Ch\^atelet surfaces (see \cite{BM1}). In \cite{BM1} we treat the case where $n\geq 9$, and in \cite{BM2} we further specialize the forms $q_i$ and prove an asymptotic for $n=8$.\\

The strategy that we adopt here builds on \cite{BM1}. There we used multiplicative characters to deal with the first equation and additive characters (the circle method) to detect the second equation. A vital `trick' was to use the modulus of the multiplicative character to reduce the size of the modulus in the circle method. This idea can also be used while applying the circle method to detect both the equations. Say we use modulus $1\leq q_1\leq B$ to detect the first equation $Q_1(\mathbf{m})=0$, which has `size' $B^2$. Then we split the second equation $Q_2(\mathbf{m})=0$ into a congruence $Q_2(\mathbf{m})\equiv 0\bmod{q_1}$ and an (integral) equation $Q_2(\mathbf{m})/q_1=0$. Now to detect the last equation by the circle method we need modulus of size $B/\sqrt{q_1}$. Hence the total modulus $q_1q_2$ has size $B^{3/2}$, instead of $B^2$ which should be the size if one used the circle method independently for both the equations. Since the size of the modulus is much smaller than the square of the length of the variables $m_i$, we save by applying Poisson summation formula to each variable. This is already sufficient to give us an asymptotic for sufficiently many variables  $n\geq 15$. But the method allows us to have a Kloosterman refinement in the first application of the circle method and a double Kloosterman refinement in the second application. This together with subconvexity for Dirichlet $L$-function reduces the number of variables to $n\geq 11$. \\

If the forms share an eigenspace then our analysis actually yields an asymptotic for $n\geq 10$. Further if we assume that both the forms are diagonal then we can have double Kloosterman refinement in both the applications of the circle method and thereby reduce the number of variables further, and get a result as strong as that of Cook \cite{C}. Finally we should mention that the nesting process, of course, comes with some extra complications beyond those reported and carefully tackled by Heath-Brown in \cite{H}. But these, as we will see, can be handled with some extra work. Also we note that for $\Delta=2$ we can use our method to improve the corresponding result of Birch. This requires a different treatment for the estimation of the character sums corresponding to `bad primes' in Section~\ref{char-sum-section}. The reader will observe that from our current estimates we get an asymptotic for $n\geq 13$ in the case $\Delta=2$. It should be possible to improve this further. \\

Theorem~\ref{mthm} gives an analytic proof of the Hasse principle for intersections of quadrics under the mild assumptions that we made. The Hasse principle is of course known in this case from the work of Colliot-Th\'el\`ene, Sansuc and Swinnerton-Dyer \cite{CSS}, \cite{CSS2}. The existence of the $p$-adic local solutions follows from the work of Demyanov \cite{De}.\\ 

Though we chose to state our result with a smooth weight function, it will be quite evident from our analysis that the smooth weight can be removed. In fact this is the reason why we have kept an extra parameter $H$ in \eqref{new-birch}. In Section~\ref{sharp-cuts} we will prove the following.

\begin{theorem}
\label{thm2}
Suppose $\mathcal P=\prod_{i=1}^n[c_i,d_i]$ is a box in $\mathbb R^n$, such that 
$\mathcal{P}\cap V^\star=\emptyset$. Then we have  
$$
\mathop{\sum\dots\sum}_{\substack{\mathbf{m}\in B\mathcal{P}\\Q_1(\mathbf{m})=Q_2(\mathbf{m})=0}}1=\mathfrak{S}J_0(\mathcal P)B^{n-4}+O\left(B^{n-4-\delta}\right)
$$
with some $\delta>0$, as long as $n\geq 11$. Here the singular integral is given by
\begin{align*}
J_0(\mathcal{P})=\mathop{\lim\lim}_{\varepsilon_1,\varepsilon_2\rightarrow 0}\frac{1}{4\varepsilon_1\varepsilon_2}\mathop{\int}_{\substack{|Q_1(\mathbf{y})|<\varepsilon_1\\|Q_2(\mathbf{y})|<\varepsilon_2\\\mathbf{y}\in\mathcal{P}}}\mathrm{d}\mathbf{y}.
\end{align*}\\
\end{theorem}

For other investigations related to systems of quadratic forms we refer the reader to the works of Dietmann \cite{D}, Heath-Brown \cite{HB3}, Schmidt \cite{S} among others. For intersections of a quadric and a cubic the method of this paper yields an asymptotic formula for the count of rational points as long as $n\geq 28$. We wish to take this up in detail in a separate paper. Unfortunately the method does not readily yield any positive result in the case of intersection of two cubics. For intersections of pairs of diagonal cubics we refer the reader to the works of Br\"udern and Wooley. 
\\

\ack The author wishes to thank Professor T.D. Browning and the anonymous referee for several helpful comments. The author is supported by SwarnaJayanti Fellowship 2011-12, DST, Government of India.\\


\section{Nested $\delta$-method and Poisson summation}

We shall now make the above argument precise. First let us briefly recall a version of the circle method introduced in \cite{DFI} and \cite{H}. The starting point is a smooth approximation of the $\delta$-symbol - $\delta(0)=1$ and $\delta(m)=0$ for $m\in \mathbb{Z}-\{0\}$. \\

\begin{lemma}
For any $Q>1$ there is a positive constant $c_Q$, and a smooth function $h(x,y)$ defined on $(0,\infty)\times\mathbb R$, such that
\begin{equation}
\delta(n)=\frac{c_Q}{Q^2}\sum_{q=1}^{\infty}\;\sideset{}{^{\star}}\sum_{a \bmod{q}}e_q\left(an\right)h\left(\frac{q}{Q},\frac{n} {Q^2}\right)\label{cm}
\end{equation}
for $n\in\mathbb Z$. Here $e_q(x)=e^{2\pi i x/q}$, the $\star$ over the sum indicates that $a$ and $q$ are coprime. The constant $c_Q$ satisfies $c_Q=1+O_A(Q^{-A})$ for any $A>0$. \\
\end{lemma}

The definition of the smooth function $h$ will be recalled in Section~\ref{integral-section}. At present we just note that $h(x,y)\ll x^{-1}$ for all $y$, and $h(x,y)$ is non-zero only for $x\leq\max\{1,2|y|\}$. So the $q$-sum in \eqref{cm} is in fact finite. In practice, to detect the equation $n=0$ for a sequence of integers in the range $[-X,X]$, it is logical to choose $Q=X^{1/2}$. \\

Let $N(B)$ denote the left hand side of \eqref{new-birch}. Using \eqref{cm} and choosing $Q$ to be $B$, we get
$$
N(B)=\frac{c_B}{B^2}\sum_{q_1=1}^{\infty}\;\sideset{}{^{\star}}\sum_{a_1 \bmod{q_1}}\mathop{\sum\dots\sum}_{\substack{\mathbf{m}\in\mathbb{Z}^n\\Q_2(\mathbf{m})=0}}
e_{q_1}\left(a_1Q_1(\mathbf{m})\right)h\left(\frac{q_1}{B},\frac{Q_1(\mathbf{m})}{B^2}\right)W\left(\frac{\mathbf{m}}{B}\right),
$$
which we rewrite as 
$$
\frac{c_B}{B^2}\sum_{q_1=1}^{\infty}\;\sideset{}{^{\star}}\sum_{a_1 \bmod{q_1}}\mathop{\sum\dots\sum}_{\substack{\mathbf{m}\in\mathbb{Z}^n\\Q_2(\mathbf{m})\equiv 0\bmod{q_1}}}e_{q_1}\left(a_1Q_1(\mathbf{m})\right)h\left(\frac{q_1}{B},\frac{Q_1(\mathbf{m})}{B^2}\right)W\left(\frac{\mathbf{m}}{B}\right)\delta\left(\frac{Q_2(\mathbf{m})}{q_1}\right).
$$
This seemingly trivial step is most vital in this paper. We are building up the second application of the circle method using the existing modulus from the first application of the circle method. This `nesting' process acts like a conductor lowering mechanism and we end up having sums over $\mathbf{m}\in \mathbb Z^n$ with modulus of size $B^{3/2}$ instead of $B^2$, which is what one has without nesting. \\

There is one technical problem that we need to take care of. For small $q_1$ the congruence condition is not strong enough to lower the conductor sufficiently. But in this case the oscillation in the analytic function $h(q_1/B,Q_1(\mathbf{m})/B^2)$ is large. So without worsening the situation we can push in a weight function that will regulate the size of $Q_2(\mathbf{m})/q_1$. Let $U$ be a bump function with support $[-2,2]$ and such that $U(x)=1$ for $x\in [-1,1]$, $U(-x)=U(x)$ and $U^{(j)}\ll_j 1$. We have
\begin{align}
\label{intro-u}
N(B)=\frac{c_B}{B^2}\sum_{q_1=1}^{\infty}\;\sideset{}{^{\star}}\sum_{a_1 \bmod{q_1}}&\mathop{\sum\dots\sum}_{\substack{\mathbf{m}\in\mathbb{Z}^n\\Q_2(\mathbf{m})\equiv 0\bmod{q_1}}}e_{q_1}\left(a_1Q_1(\mathbf{m})\right)h\left(\frac{q_1}{B},\frac{Q_1(\mathbf{m})}{B^2}\right)\\
\nonumber &\times W\left(\frac{\mathbf{m}}{B}\right)\delta\left(\frac{Q_2(\mathbf{m})}{q_1}\right)\:U\left(\frac{Q_2(\mathbf{m})}{q_1B}\right).
\end{align}
Applying \eqref{cm} again, choosing 
$
Q=\sqrt{B}
$
in this application of the circle method, and rearranging the sums, we get
\begin{align*}
N(B)=\frac{1}{B^3}\sum_{q_1=1}^{\infty}\;\sum_{q_2=1}^{\infty}\;\sideset{}{^{\star}}\sum_{a_1 \bmod{q_1}}\;\sideset{}{^{\star}}\sum_{a_2 \bmod{q_2}}\;&N(\mathbf{a},\mathbf{q};B)+O(B^{-2013})
\end{align*}
where
\begin{align*}
N(\mathbf{a},\mathbf{q};B)=&\mathop{\sum\dots\sum}_{\substack{\mathbf{m}\in\mathbb{Z}^n\\Q_2(\mathbf{m})\equiv 0\bmod{q_1}}}e_{q_1q_2}\left(a_1Q_1(\mathbf{m})q_2+a_2Q_2(\mathbf{m})\right)\\
&\times h\left(\frac{q_1}{B},\frac{Q_1(\mathbf{m})}{B^2}\right)h\left(\frac{q_2}{\sqrt{B}},\frac{Q_2(\mathbf{m})}{q_1B}\right)U\left(\frac{Q_2(\mathbf{m})}{q_1B}\right)W\left(\frac{\mathbf{m}}{B}\right).
\end{align*}
Applying the Poisson summation formula with modulus $q_1q_2$ we deduce that
\begin{align}
\label{poisson}
N(\mathbf{a},\mathbf{q};B)=&\frac{B^n}{(q_1q_2)^n}\mathop{\sum\dots\sum}_{\mathbf{m}\in\mathbb{Z}^n} I_{q_1,q_2}(\mathbf{m})\\
\nonumber \times & \left[\mathop{\sum\dots\sum}_{\substack{\mathbf{b}\bmod{q_1q_2}\\Q_2(\mathbf{b})\equiv 0\bmod{q_1}}}e_{q_1q_2}\left(a_1Q_1(\mathbf{b})q_2+a_2Q_2(\mathbf{b})+\mathbf{b}.\mathbf{m}\right)\right],
\end{align}
where
\begin{align}
\label{integral}
I_{q_1,q_2}(\mathbf{m})=&\int_{\mathbb R^n}h\left(\frac{q_1}{B},Q_1(\mathbf{y})\right)h\left(\frac{q_2}{\sqrt{B}},\frac{BQ_2(\mathbf{y})}{q_1}\right)\\
\nonumber &\times U\left(\frac{BQ_2(\mathbf{y})}{q_1}\right)W\left(\mathbf{y}\right)e_{q_1q_2}\left(-B\mathbf{m}.\mathbf{y}\right)\mathrm{d}\mathbf{y}.
\end{align}
We set
\begin{align}
\label{char-sum}
C_{q_1,q_2}(\mathbf{m})=\sideset{}{^{\star}}\sum_{a_1 \bmod{q_1}}\;\sideset{}{^{\star}}\sum_{a_2 \bmod{q_2}}\mathop{\sum\dots\sum}_{\substack{\mathbf{b}\bmod{q_1q_2}\\Q_2(\mathbf{b})\equiv 0\bmod{q_1}}}e_{q_1q_2}\left(a_1Q_1(\mathbf{b})q_2+a_2Q_2(\mathbf{b})+\mathbf{b}.\mathbf{m}\right).
\end{align}\\

\begin{lemma}
\label{lemma2}
We have
$$
N(B)=N^\star(B)+O(B^{-2013})
$$
where
\begin{align*}
N^\star(B)=B^{n-3}\mathop{\sum\sum}_{\substack{q_1\ll B\\q_2\ll B^{1/2}}}\;\frac{1}{(q_1q_2)^n}&\mathop{\sum\dots\sum}_{\mathbf{m}\in\mathbb{Z}^n} \:C_{q_1,q_2}(\mathbf{m})I_{q_1,q_2}(\mathbf{m}).
\end{align*}\\
\end{lemma}

To complete the proof of the lemma all that remains to show is that the contribution from $q_1\gg B$ and $q_2\gg B^{1/2}$ is negligibly small for some appropriate implied constants. Indeed 
$$
h\left(\frac{q_1}{B},Q_1(\mathbf{y})\right)
$$ 
vanishes for 
$$
q_1>c_1B\;\;\;\text{where}\;\;\;c_1=1+\sup_{\mathbf{y}\in \text{Supp}(W)}| Q_1(\mathbf{y})|,
$$
and
$$
h\left(\frac{q_2}{\sqrt{B}},\frac{BQ_2(\mathbf{y})}{q_1}\right)U\left(\frac{BQ_2(\mathbf{y})}{q_1}\right)
$$
vanishes for $q_2>2B^{1/2}$. Note that without the weight function $U$, we would need to take $q_2$ upto $B^{3/2}/q_1$, which can be as large as $B^{3/2}$ for small $q_1$.\\


\section{The character sum}
\label{char-sum-section}

The character sum \eqref{char-sum} can be compared with the character sum which appears in Lemma~8 of \cite{BM1}. So we expect to have analogous results for this character sum. In particular we have multiplicativity property. Let $(q_1,q_2)=d$, and we write $q_1=q_1'd_1$ and $q_2=q_2'd_2$ where $d_1d_2|d^\infty$ and $(q_1'q_2',d)=1$ (consequently $(q_1',q_2')=1$). Then we have
$$
C_{q_1,q_2}(\mathbf{m})=C_{q_1',1}(\mathbf{m})C_{1,q_2'}(\mathbf{m})C_{d_1,d_2}(\mathbf{m}).
$$\\

A major part of \cite{BM1} is devoted to the study of the character sum
$$
S_{q_1,q_2}(\mathbf{m})=\sideset{}{^{\star}}\sum_{a \bmod{q_2}}\mathop{\sum\dots\sum}_{\substack{\mathbf{b}\bmod{q_1q_2}\\Q_2(\mathbf{b})\equiv 0\bmod{q_1}\\Q_1(\mathbf{b})\equiv 0\bmod{q_2}}}e_{q_1q_2}\left(aQ_2(\mathbf{b})+\mathbf{b}.\mathbf{m}\right),
$$ 
where $q_1$ and $q_2$ positive integers, and $\mathbf{m}\in \mathbb{Z}^n$. (We introduce the convention that $S_{q_1,q_2}(\mathbf{m})=0$ if either $q_1$ or $q_2$ is not a positive integer or if $\mathbf{m}$ is not an integral vector.) This is similar to the character sum $C_{q_1,q_2}(\mathbf{m})$. We will analyse this character sum using the same techniques that we employed in Sections~4, 5 and 6 of \cite{BM1}. But first we will explicitly relate the two character sums, which will enable us to translate some of the end results from \cite{BM1} to the current situation. This shortens the exposition in this section tremendously, albeit at the cost of making it not entirely self-contained. We begin by proving the following relation.\\

\begin{lemma}
\label{lemma3}
Let $p$ be any prime and $r,\ell\geq 0$, then we have
\begin{align}
\label{ship1}
C_{p^r,p^\ell}(\mathbf{m})=p^rS_{p^r,p^\ell}(\mathbf{m})-p^{r-2}S_{p^{r-1},p^{\ell+1}}(\mathbf{m})-p^{n+r-2}S_{p^{r-1},p^{-\ell}}(p^{-1}\mathbf{m}).
\end{align}\\
\end{lemma}  
 
For any prime $p$ we have
\begin{align}
\label{prime-c-exp}
C_{p^r,p^\ell}(\mathbf{m})=\sideset{}{^{\star}}\sum_{a_2 \bmod{p^\ell}}\mathop{\sum\dots\sum}_{\substack{\mathbf{b}\bmod{p^{r+\ell}}\\Q_2(\mathbf{b})\equiv 0\bmod{p^r}}}c_{p^r}\left(Q_1(\mathbf{b})\right)e_{p^{r+\ell}}\left(a_2Q_2(\mathbf{b})+\mathbf{b}.\mathbf{m}\right)
\end{align}
where 
$$
c_q(m)=\sideset{}{^\star}\sum_{a\bmod{q}}e_q(am)
$$ denotes the Ramanujan sum. 
Using the well known formula for the Ramanujan sum 
$$
c_q(m)=\sum_{d|(m,q)}\mu(q/d)d
$$
(where $\mu$ denotes the M\"{o}bius function) we get
\begin{align}
\label{break-up}
p^r\sideset{}{^{\star}}\sum_{a_2 \bmod{p^\ell}}\mathop{\sum\dots\sum}_{\substack{\mathbf{b}\bmod{p^{r+\ell}}\\Q_2(\mathbf{b})\equiv 0\bmod{p^r}\\Q_1(\mathbf{b})\equiv 0\bmod{p^r}}}&e_{p^{r+\ell}}\left(a_2Q_2(\mathbf{b})+\mathbf{b}.\mathbf{m}\right)\\
\nonumber &-p^{r-1}\sideset{}{^{\star}}\sum_{a_2 \bmod{p^\ell}}\mathop{\sum\dots\sum}_{\substack{\mathbf{b}\bmod{p^{r+\ell}}\\Q_2(\mathbf{b})\equiv 0\bmod{p^r}\\Q_1(\mathbf{b})\equiv 0\bmod{p^{r-1}}}}e_{p^{r+\ell}}\left(a_2Q_2(\mathbf{b})+\mathbf{b}.\mathbf{m}\right)
\end{align}
with the understanding that the second term appears only in the case $r\geq 1$.  So the first term in \eqref{break-up} is just $p^rS_{p^r,p^\ell}(\mathbf{m})$. In the second term in \eqref{break-up} the sum over $a_2$ yields the Ramanujan sum $c_{p^\ell}(Q_2(\mathbf{b})/p^r)$. It follows that this term equals
\begin{align}
\label{relate}
p^{r+\ell-1}\mathop{\sum\dots\sum}_{\substack{\mathbf{b}\bmod{p^{r+\ell}}\\Q_2(\mathbf{b})\equiv 0\bmod{p^{r+\ell}}\\Q_1(\mathbf{b})\equiv 0\bmod{p^{r-1}}}}e_{p^{r+\ell}}\left(\mathbf{b}.\mathbf{m}\right)-p^{r+\ell-2}\mathop{\sum\dots\sum}_{\substack{\mathbf{b}\bmod{p^{r+\ell}}\\Q_2(\mathbf{b})\equiv 0\bmod{p^{r+\ell-1}}\\Q_1(\mathbf{b})\equiv 0\bmod{p^{r-1}}}}e_{p^{r+\ell}}\left(\mathbf{b}.\mathbf{m}\right).
\end{align}
Again the second sum appears only if $\ell\geq 1$. Now 
$$
S_{p^{r-1},p^{\ell+1}}(\mathbf{m})=\mathop{\sum\dots\sum}_{\substack{\mathbf{b}\bmod{p^{r+\ell}}\\Q_2(\mathbf{b})\equiv 0\bmod{p^{r-1}}\\Q_1(\mathbf{b})\equiv 0\bmod{p^{r-1}}}}c_{p^{\ell+1}}(Q_2(\mathbf{b})/p^{r-1})e_{p^{r+\ell}}\left(\mathbf{b}.\mathbf{m}\right),
$$ 
which decomposes as 
\begin{align}
\label{relate2}
p^{\ell+1}\mathop{\sum\dots\sum}_{\substack{\mathbf{b}\bmod{p^{r+\ell}}\\Q_2(\mathbf{b})\equiv 0\bmod{p^{r+\ell}}\\Q_1(\mathbf{b})\equiv 0\bmod{p^{r-1}}}}e_{p^{r+\ell}}\left(\mathbf{b}.\mathbf{m}\right)-p^{\ell}\mathop{\sum\dots\sum}_{\substack{\mathbf{b}\bmod{p^{r+\ell}}\\Q_2(\mathbf{b})\equiv 0\bmod{p^{r+\ell-1}}\\Q_1(\mathbf{b})\equiv 0\bmod{p^{r-1}}}}e_{p^{r+\ell}}\left(\mathbf{b}.\mathbf{m}\right).
\end{align}
Comparing this with \eqref{relate} we conclude that for $\ell\geq 1$ we have
\begin{align*}
C_{p^r,p^\ell}(\mathbf{m})=p^rS_{p^r,p^\ell}(\mathbf{m})-p^{r-2}S_{p^{r-1},p^{\ell+1}}(\mathbf{m}).
\end{align*}
In particular for $r=0$ we have $C_{1,p^\ell}(\mathbf{m})=S_{1,p^\ell}(\mathbf{m})$. For $\ell=0$ the second term in \eqref{relate2} vanishes unless $p|\mathbf{m}$ in which case it is given by $p^nS_{p^{r-1},1}(p^{-1}\mathbf{m})$. So it follows that
$$
C_{p^r,1}(\mathbf{m})=p^rS_{p^r,1}(\mathbf{m})-p^{r-2}S_{p^{r-1},p}(\mathbf{m})-p^{n+r-2}S_{p^{r-1},1}(p^{-1}\mathbf{m}).
$$ 
This concludes the proof of Lemma~\ref{lemma3}.\\

We can now freely borrow the results from \cite{BM1} about the character sum $S_{q_1,q_2}(\mathbf{m})$ to obtain sufficient bounds for the sum $C_{q_1,q_2}(\mathbf{m})$. Let $M_i$ be the symmetric matrices associated to the quadratic forms $Q_i$. Let $Q_2^\star$ denote the adjoint quadratic form with underlying matrix $M_2^\star=(\det M_2) M_2^{-1}$. (Recall that $M_2$ is invertible.) For any odd prime $p$ let $\chi_p$ be the unique quadratic character modulo $p$, and let $\varepsilon(p)=1$ for $p\equiv 1\bmod{4}$ and $\varepsilon(p)=i$ for $p\equiv 3\bmod{4}$. Also let $g_q(m)$ denote the Gauss sum with modulus $q$. Lemma~15 and (4.1) of \cite{BM1} translate to yield the following, as from Lemma~\ref{lemma3} we have $C_{1,q_2}(\mathbf{m})=S_{1,q_2}(\mathbf{m})$. \\

\begin{lemma}
\label{lemma4}
Let $p$ be prime with $p\nmid 2\det M_2$. For even $n$ we have
\begin{align*}
C_{1,p^\ell}(\mathbf{m})=\varepsilon(p)^{n\ell}\chi_p(\det M_2)^\ell p^{n\ell/2}c_{p^\ell}(Q_2^\star(\mathbf{m})).
\end{align*}
For odd $n$ we have
\begin{align}
\label{char-sum-bd-00}
C_{1,p^\ell}(\mathbf{m})=
\begin{cases}p^{n\ell/2}c_{p^\ell}(Q_2^\star(\mathbf{m})) &\text{if $\ell$ is even}\\
\varepsilon(p)^n\chi_p(-1) p^{n\ell/2}g_{p^\ell}(Q_2^\star(\mathbf{m})) &\text{if $\ell$ is odd}.\end{cases}
\end{align}
Finally for any prime $p$ and any integral vector $\mathbf{m}$ we have
$$
C_{1,p^\ell}(\mathbf{m})\ll p^{\ell(n/2+1)}.
$$\\
\end{lemma}

In Section~4 of \cite{BM1} we also showed that there is cancellation when we average $C_{1,q}(\mathbf{m})$ over $q$. This follows from the oscillation in the sign of the Gauss sum, and can be shown using the standard properties of the Dirichlet L-function. One should recall that hybrid subconvexity of the Dirichlet $L$-function plays a crucial role in this deduction. To state the result we define 
\begin{align}
\label{nm}
N=N(\mathbf{m})=\begin{cases} 2Q_2^\star(\mathbf{m})\det M_2 &\text{if $Q_2^\star(\mathbf{m})\neq 0$};\\
2\det M_2 &\text{otherwise}.
\end{cases}
\end{align}
The following lemma summarizes the content of Section~4 of \cite{BM1}. \\

\begin{lemma}
\label{lemma5} 
Let $M\in\mathbb N$ with $N|M$. Suppose $n$ is odd then we have 
\begin{align}
\label{char-sum-bd-1}
\sum_{\substack{q_2\leq x\\(q_2,M)=1}}C_{1,q_2}(\mathbf{m})\ll
\begin{cases}|\mathbf{m}|^{7/16+\varepsilon}x^{(n+2)/2+\varepsilon}M^\varepsilon &\text{if $(-1)^{(n-1)/2}Q_2^\star(\mathbf{m})\neq \Box$}\\
x^{(n+3)/2+\varepsilon}M^\varepsilon &\text{if $(-1)^{(n-1)/2}Q_2^\star(\mathbf{m})=\Box$}.\end{cases}
\end{align}
Suppose $n$ is even then we have
\begin{align}
\label{char-sum-bd-0}
\sum_{\substack{q_2\leq x\\(q_2,M)=1}}C_{1,q_2}(\mathbf{m})\ll
\begin{cases}x^{(n+2)/2+\varepsilon}M^\varepsilon &\text{if $Q_2^\star(\mathbf{m})\neq 0$}\\
x^{(n+3)/2+\varepsilon}M^\varepsilon &\text{if $Q_2^\star(\mathbf{m})=0$ and $(-1)^{n/2}\det M_2\neq \Box$}\\
x^{(n+4)/2+\varepsilon}M^\varepsilon &\text{if $Q_2^\star(\mathbf{m})=0$ and $(-1)^{n/2}\det M_2=\Box$}.\end{cases}
\end{align}
\\
\end{lemma}

We now turn to derive a bound for the sum $C_{q_1,1}(\mathbf{m})$. Lemma~22 of \cite{BM1} gives 
$$C_{q_1,1}(\mathbf{m})\ll (q_1,\Delta_V^\infty)^{n/2-2}(q_1,\mathbf{m})^{n/2-2}q_1^{n/2+1+\varepsilon}.$$ But here we will be able to remove the extra factor $(q_1,\Delta_V^\infty)^{n/2-2}(q_1,\mathbf{m})^{n/2-2}$ from the bound. We begin by recalling the definition of $\Delta_V$. Let $V$ be the non-singular projective variety defined by $Q_1(\mathbf{m})=Q_2(\mathbf{m})=0$ in the space $\mathbb{P}^{n-1}$. The dual variety $V^\star$ is given by an equation $G(\mathbf{m})=0$ where $G$ is an irreducible form with integral coefficients of degree $4(n-2)$ (see Section~2.2 of \cite{BM1}). For any point $\mathbf{b}=[b_1,b_2]\in\mathbb{P}^1_\mathbb{Q}$ consider the matrix $M(\mathbf{b})=b_1M_1+b_2M_2$. From Lemma~1.13 of \cite{CSS} it follows that $\text{rank}\:M(\mathbf{b})\geq n-1$ for all $\mathbf{b}\in\mathbb{P}^1$. This property (as well as the non-singularity of $V$) continues to hold modulo $p$ for all but finitely many primes $p$. Let $\mathcal{B}$ be the product of this finite set of primes. Furthermore, from Ried's thesis \cite{R} it follows that the binary form $\det M(\mathbf{b})$ has non-zero discriminant $\mathcal{D}$. Let 
\begin{align}
\label{deltav}
\Delta_V=2(\det M_2)\:\mathcal{B}\:\mathcal{D}.
\end{align}
By `bad primes' we shall mean the primes dividing $\Delta_V$.\\

\begin{lemma}
\label{lemma6}
We have
\begin{align}
\label{char-sum-bd-3b}
C_{q_1,1}(\mathbf{m})\ll q_1^{n/2+1+\varepsilon},
\end{align}
where the implied constant depends on $Q_i$ and $\varepsilon$.\\
\end{lemma}

Indeed we have
\begin{align}
\label{middd}
C_{p^r,1}(\mathbf{m})=p^{-r}\sideset{}{^{\star}}\sum_{a_1 \bmod{p^r}}\sum_{c\bmod{p^r}}\mathop{\sum\dots\sum}_{\substack{\mathbf{b}\bmod{p^{r}}}}
e_{p^{r}}\left(a_1Q_1(\mathbf{b})+a_1cQ_2(\mathbf{b})+\mathbf{b}.\mathbf{m}\right).
\end{align}
To which we apply Cauchy's inequality and Lemma~13 of \cite{BM1} to get
\begin{align*}
C_{p^r,1}(\mathbf{m})\leq p^{nr/2-r}\sideset{}{^{\star}}\sum_{a_1 \bmod{p^r}}\sum_{c\bmod{p^r}} K_{p^r}(c)^{1/2}
\end{align*}
where 
$$
K_{p^r}(c)=1+\#\{\mathbf{v}\bmod{p^r}:M(c)\mathbf{v}\equiv\mathbf{0}\bmod{p^r}\},
$$
with $M(c)=M_1+cM_2$. This can be compared with (5.3) of \cite{BM1}. Notice that here we have $h=1$ as from the beginning $a_1$ is coprime with $p$. Using Lemma~21 of \cite{BM1} we get \eqref{char-sum-bd-3b} if $(q_1,\Delta_V)=1$.\\

Now take $p|\Delta_V$, $r\geq 1$, and consider the inner sum over $\mathbf{b}$ in \eqref{middd}. Taking $j=[(r+1)/2]$ and writing $\mathbf{b}=\mathbf{b}_1+\mathbf{b}_2p^j$ we see that the sum over $\mathbf{b}$ is given by
$$
\mathop{\sum\dots\sum}_{\substack{\mathbf{b}_1\bmod{p^{j}}}}
e_{p^{r}}\left(a_1~^t\mathbf{b}_1M(c)\mathbf{b}_1+\mathbf{b}_1.\mathbf{m}\right)\mathop{\sum\dots\sum}_{\substack{\mathbf{b}_2\bmod{p^{r-j}}}}
e_{p^{r-j}}\left(2a_1~^t\mathbf{b}_1M(c)\mathbf{b}_2+\mathbf{b}_2.\mathbf{m}\right).
$$
This reduces to
$$
p^{n(r-j)}\mathop{\sum\dots\sum}_{\substack{\mathbf{b}_1\bmod{p^{j}}\\2a_1M(c)\mathbf{b}_1+\mathbf{m}\equiv \mathbf{0}\bmod{p^{r-j}}}}
e_{p^{r}}\left(a_1~^t\mathbf{b}_1M(c)\mathbf{b}_1+\mathbf{b}_1.\mathbf{m}\right),
$$
which is dominated by
$$
p^{n(r-j)}\mathop{\sum\dots\sum}_{\substack{\mathbf{b}_1\bmod{p^{j}}\\2a_1M(c)\mathbf{b}_1+\mathbf{m}\equiv \mathbf{0}\bmod{p^{r-j}}}}1\ll_V p^{nr/2}\:K_{p^j}(c)
$$
Consequently we have
$$
C_{p^r,1}(\mathbf{m})\ll p^{nr/2}\sum_{c\bmod{p^r}}K_{p^j}(c).
$$
Let $K$ be a finite extension of the local field $\mathbb{Q}_p$ containing all the eigenvalues of the symmetric integral matrix $M_2^\star M_1$.  Let $\pi$ be an uniformizer for this extension, so that $p\mathcal{O}_K=\pi^e\mathcal{O}_K$ where $e$ is the ramification index. Now we have an invertible matrix $U$ and a diagonal matrix  $D=\text{diag}(d_1,\dots,d_n)$ with entries in $\mathcal{O}_K$, such that $~^tUU=\pi^hI$ for some non-negative integer $h$ and such that $\pi^{h}M_2^\star M_1=~^tUDU$. The entries $d_i$ of the diagonal matrix $D$ are in fact the eigenvalues of the matrix $M_2^\star M_1$, with associated eigenvector $~^tU\mathbf{e}_i$ (where $\mathbf{e}_i$ is the $i$-th standard basis vector). Since we are assuming that the singular locus has dimension $\Delta=1$, it follows that the eigenvalues  $d_i$ are distinct. Now $M(c)\mathbf{v}\equiv\mathbf{0}\bmod{p^j}$ implies $M_2^\star M_1 \mathbf{v}\equiv -c(\det M_2)\mathbf{v}\bmod{p^j}$, which in turn implies 
$$
DU \mathbf{v}\equiv -c(\det M_2)U\mathbf{v}\bmod{\pi^{ej-h}}.
$$
Writing $U\mathbf{v}=~^t(u_1,\dots,u_n)$ we get a collection of independent congruences 
$$
d_iu_i\equiv -c(\det M_2)u_i \bmod{\pi^{ej-h}}.
$$ 
This leads us to consider the $p$-adic distance $|d_i+c\det M_2|_p$. Suppose $\delta=\min_{d_i\neq d_j}|d_i-d_j|_p$. If $|d_i+c\det M_2|_p<\delta$ then by Krasner's lemma we have $d_i\in \mathbb{Z}_p$. We conclude that if $i$ is such that $d_i\notin \mathbb{Z}_p$ then the number of $u_i$ satisfying the above congruence is $O_V(1)$ (where the implied constant does not depend on $c$ or $j$). On the other hand if $d_i\in\mathbb{Z}_p$ then we can suppose that $~^tU\mathbf{e}_i\in \mathbb{Z}_p^n$, and consequently $u_i\in \mathbb{Z}_p$. Now for any $c$ the number of $u_i\in \mathbb{Z}/p^{j}\mathbb{Z}$ satisfying the $i$-th congruence is bounded by $O_V((d_i+c\det M_2,p^j))$. From $\mathbf{u}$ we obtain $\mathbf{v}$ by multiplying with $~^tU$. Hence we are able to conclude that
$$
\sum_{c\bmod{p^r}}K_{p^j}(c)\ll \sum_{c\bmod{p^r}} \prod_{d_i\in\mathbb{Z}_p}(d_i+c\det M_2,p^j).
$$
Now using the fact that the multiplicities of the eigenvalues are $1$ we derive that
$$
\sum_{c\bmod{p^r}}K_{p^j}(c)\ll \sum_{d_i\in\mathbb{Z}_p}\sum_{c\bmod{p^r}}(d_i+c\det M_2,p^j).
$$
We have
$$
\sum_{c\bmod{p^r}}(d_i+c\det M_2,p^j)\ll \sum_{a=0}^jp^{a}\sum_{\substack{c\bmod{p^r}\\c\det M_2\equiv -d_i\bmod{p^a}}}1\ll p^{r}.
$$
As result we conclude that \eqref{char-sum-bd-3b} holds for any $p^r|\Delta_V^\infty$. This completes the proof of Lemma~\ref{lemma6}.\\

Next we will obtain a bound for $C_{q_1,1}(\mathbf{m})$ on average over $q_1$. We will prove the following analogue of Lemmas~23 and 24 of \cite{BM1}.\\

\begin{lemma}
\label{lemma7}
For $G(\mathbf{m})\neq 0$ we have
\begin{align}
\label{char-sum-bd-2}
\sum_{\substack{q_1\leq x\\(q_1,\Delta_VG(\mathbf{m}))=1}}\left|C_{q_1,1}(\mathbf{m})\right|\ll x^{(n+2)/2+\varepsilon}.
\end{align}
For $G(\mathbf{m})=0$ we have
\begin{align}
\label{char-sum-bd-3}
\sum_{\substack{q_1\leq x\\(q_1,\Delta_V\mathbf{m})=1}}\left|C_{q_1,1}(\mathbf{m})\right|\ll x^{(n+3)/2+\varepsilon}.
\end{align}\\
\end{lemma}

First consider the case where $p\nmid \Delta_VG(\mathbf{m})$ (so in particular $G(\mathbf{m})\neq 0$ and $p\nmid \mathbf{m}$). Then from Lemma~\ref{lemma3} we get
$$
C_{p^r,1}(\mathbf{m})=p^rS_{p^r,1}(\mathbf{m})-p^{r-2}S_{p^{r-1},p}(\mathbf{m}).
$$
From Lemmas~19 and 20 of \cite{BM1} we get $p^rS_{p^r,1}(\mathbf{m})=0$ if $r>1$, and $pS_{p,1}(\mathbf{m})\ll p^{n/2}$.  Also from Lemma~\ref{lemma4} above and Lemma~27 of \cite{BM1} we obtain $ p^{r-2}S_{p^{r-1},p}(\mathbf{m})\ll p^{-1}p^{r(n/2+1)}$. Consequently we get
$$
C_{q_1,1}(\mathbf{m})\ll \frac{q_1^{n/2+1+\varepsilon}}{\text{rad}(q_1)},
$$
where $\text{rad}(m)$ is the largest square-free factor of the positive integer $m$. From this we derive \eqref{char-sum-bd-2}. \\

Now to prove \eqref{char-sum-bd-3} we note that from Lemmas~19 of \cite{BM1} we have $pS_{p,1}(\mathbf{m})\ll p^{(n+1)/2}$. This together with Lemma~\ref{lemma4} yield $C_{p,1}(\mathbf{m})\ll p^{(n+1)/2}$. For $r\geq 2$ we use the bound from Lemma~\ref{lemma6} above. We get
\begin{align*}
\sum_{\substack{q_1\leq x\\(q_1,\Delta_V\mathbf{m})=1}}\left|C_{q_1,1}(\mathbf{m})\right|&\ll \sum_{\substack{q_1\leq x\\(q_1,\Delta_V\mathbf{m})=1\\q_1\:\Box-\text{free}}}\left|C_{q_1,1}(\mathbf{m})\right|\sum_{\substack{q_1'\leq x/q_1\\(q_1',\Delta_V\mathbf{m})=1\\q_1'\:\Box-\text{full}}}\left|C_{q_1',1}(\mathbf{m})\right|\\
&\ll x^\varepsilon\sum_{\substack{q_1\leq x\\(q_1,\Delta_V\mathbf{m})=1\\q_1\:\Box-\text{free}}}\left|C_{q_1,1}(\mathbf{m})\right|\:\frac{x^{n/2+1}}{q_1^{n/2+1}}\:\frac{x^{1/2}}{q_1^{1/2}}\ll x^{(n+3)/2+\varepsilon}.
\end{align*}
This completes the proof of Lemma~\ref{lemma7}.\\


Now we consider the mixed character sum $C_{p^r,p^\ell}(\mathbf{m})$ with $r,\ell\geq 1$. We will establish the following bound.\\

\begin{lemma}
\label{lemma8}
For $d_1$, $d_2$ positive integers with $d_1|d_2^\infty$ and $d_2|d_1^\infty$, we have
\begin{align}
\label{char-sum-bd-4}
C_{d_1,d_2}(\mathbf{m})\ll (d_1d_2)^{n/2+1+\varepsilon}\frac{\left(\mathrm{rad}(d_2),Q_2^\star(\mathbf{m})\right)}{\mathrm{rad}(d_2)}.
\end{align}\\
\end{lemma}

For $p\nmid 2\det M_2$, using Lemmas~26 and 27 of \cite{BM1} we get that the first term of \eqref{ship1} is bounded as
$$
p^rS_{p^r,p^\ell}(\mathbf{m})\ll p^{(r+\ell)(n/2+1)}\frac{(p,Q_2^\star(\mathbf{m}))}{p}.
$$
The same bound also holds for the second term of \eqref{ship1}, because of same reasons if $r\geq 2$, and because of \eqref{char-sum-bd-00} in the case $r=1$. Also from Lemma~26 of \cite{BM1} we see that the bound remains valid even if $p|2\det M_2$ in the case $\ell>r$. Hence 
$$
C_{p^r,p^\ell}(\mathbf{m})\ll p^{(r+\ell)(n/2+1)}\frac{(p,Q_2^\star(\mathbf{m}))}{p}.
$$
if $p\nmid 2\det M_2$ or if $\ell>r$. Now take a prime $p|2\det M_2$ and take $r\geq \ell$. Writing $\mathbf{b}=\mathbf{u}+\mathbf{v}p^r$ in \eqref{prime-c-exp} we get 
\begin{align*}
C_{p^r,p^\ell}(\mathbf{m})=\sideset{}{^{\star}}\sum_{a_2 \bmod{p^\ell}}&\mathop{\sum\dots\sum}_{\substack{\mathbf{u}\bmod{p^{r}}\\Q_2(\mathbf{u})\equiv 0\bmod{p^r}}}c_{p^r}\left(Q_1(\mathbf{u})\right)e_{p^{r+\ell}}\left(a_2Q_2(\mathbf{u})+\mathbf{u}.\mathbf{m}\right)\\
&\times \mathop{\sum\dots\sum}_{\substack{\mathbf{v}\bmod{p^{\ell}}}}e_{p^{\ell}}\left(2a_2\nabla Q_2(\mathbf{u}).\mathbf{v}+\mathbf{v}.\mathbf{m}\right)
\end{align*}
which reduces to
\begin{align*}
C_{p^r,p^\ell}(\mathbf{m})=p^{n\ell}\sideset{}{^{\star}}\sum_{a_2 \bmod{p^\ell}}&\mathop{\sum\dots\sum}_{\substack{\mathbf{u}\bmod{p^{r}}\\Q_2(\mathbf{u})\equiv 0\bmod{p^r}\\2a_2M_2\mathbf{u}+\mathbf{m}\equiv \mathbf{0}\bmod{p^\ell}}}c_{p^r}\left(Q_1(\mathbf{u})\right)e_{p^{r+\ell}}\left(a_2Q_2(\mathbf{u})+\mathbf{u}.\mathbf{m}\right).
\end{align*}
Detecting the congruence modulo $p^r$ using exponential sum we get
\begin{align*}
p^{n\ell-r}\mathop{\sideset{}{^{\star}}\sum\sideset{}{^{\star}}\sum}_{\substack{a_1\bmod{p^r}\\a_2 \bmod{p^{r+\ell}}}}&\mathop{\sum\dots\sum}_{\substack{\mathbf{u}\bmod{p^{r}}\\2a_2M_2\mathbf{u}+\mathbf{m}\equiv \mathbf{0}\bmod{p^\ell}}}e_{p^{r+\ell}}\left(a_1p^\ell Q_1(\mathbf{u})+a_2Q_2(\mathbf{u})+\mathbf{u}.\mathbf{m}\right).
\end{align*}
The inner sum over $\mathbf{u}$ is dominated by
\begin{align*}
\mathop{\sum\dots\sum}_{\substack{\mathbf{w}\bmod{p^{\ell}}\\2a_2M_2\mathbf{w}+\mathbf{m}\equiv \mathbf{0}\bmod{p^\ell}}}\Bigl|\mathop{\sum\dots\sum}_{\substack{\mathbf{u}\bmod{p^{r-\ell}}}}e_{p^{r-\ell}}\left(~^t\mathbf{u}M(a_1,a_2)\mathbf{u}+\mathbf{u}.\mathbf{m}'\right)\Bigr|.
\end{align*}
where $M(a_1,a_2)=a_1p^\ell M_1+a_2M_2$ and $\mathbf{m}'=2a_1M_1\mathbf{w}+p^{-\ell}(2a_2M_2\mathbf{w}+\mathbf{m})$. Using Lemma~13 of \cite{BM1} we see that this is bounded by
\begin{align*}
p^{(r-\ell)n/2}\mathop{\sum\dots\sum}_{\substack{\mathbf{w}\bmod{p^{\ell}}\\2a_2M_2\mathbf{w}+\mathbf{m}\equiv \mathbf{0}\bmod{p^\ell}}}K_{p^{r-\ell}}(2M(a_1,a_2))^{1/2}.
\end{align*}
Consequently
\begin{align*}
C_{p^r,p^\ell}(\mathbf{m})\ll p^{(r+\ell)n/2-r} (2\det M_2,p^\ell)^n\mathop{\sideset{}{^{\star}}\sum\sideset{}{^{\star}}\sum}_{\substack{a_1\bmod{p^r}\\a_2 \bmod{p^{r+\ell}}}}K_{p^{r-\ell}}(2M(a_1,a_2))^{1/2}.
\end{align*}
It remains to estimate the last sum. For any $\mathbf{v}\bmod{p^{r-\ell}}$ we write $\mathbf{v}=p^c\mathbf{v}'$ with $(p,\mathbf{v}')=1$. The equation $2a_1p^\ell M_1\mathbf{v}\equiv -2a_2M_2\mathbf{v}\bmod{p^{r-\ell}}$ implies $2a_1p^{\ell+c} M_2^\star M_1\mathbf{v}'\equiv 2a_2p^c\det M_2 \mathbf{v}'\bmod{p^{r-\ell}}$. If $p^\ell\nmid \det M_2$ then we should necessarily have $c\geq r-\ell-v_p(\det M_2)$, and it follows that $K_{p^{r-\ell}}(2M(a_1,a_2))\ll 1$. Now suppose $p^\ell|\det M_2$ (then $p^\ell\ll 1$). Then using Lemma~21 of \cite{BM1} we get that 
\begin{align*}
\mathop{\sideset{}{^{\star}}\sum\sideset{}{^{\star}}\sum}_{\substack{a_1\bmod{p^r}\\a_2 \bmod{p^{r+\ell}}}}K_{p^{r-\ell}}(2M(a_1,a_2))^{1/2}\ll p^{2r}.
\end{align*}
The lemma follows.\\ 



\section{The integral}\label{integral-section}

We will follow the treatment given in Sections 7 and 8 of \cite{H}. But there are new complications, arising from the double application of the circle method, that we need to discuss. To this end, we recall some facts about the smooth function $h(x,y)$ which appears in \eqref{cm}. We begin by noting that
\begin{align}
\label{hxy-form}
h(x,y)=\sum_{j=1}^\infty \frac{1}{xj}\left\{w(xj)-w(|y|/xj)\right\}
\end{align}
where $w$ is a smooth function supported on $[1/2,1]$ with $\int w(x)\mathrm{d}x=1$. So $h(x,y)=0$ if $x\geq \max\{1,2|y|\}$. It follows that $h(x,y)$ satisfies (see \cite{H}) 
\begin{eqnarray}
x^{i} \frac{\partial^i}{\partial x^i}h(x,y)\ll_i x^{-1} & \textnormal{and} & \frac{\partial}{\partial y}h(x,y)=0 \label{hbound1}
\end{eqnarray}
for $x\leq 1$ and $|y|\leq x/2$. For $|y|>x/2$, we have
\begin{equation}
\frac{\partial^{i}}{\partial x^i}\frac{\partial^{j}}{\partial y^j}h(x,y)\ll_{i,j} x^{-1-i}|y|^{-j}. \label{hbound2}
\end{equation}
Also for any $(x,y)$, and non-negative integers $i$, $j$ and $N$, we have
\begin{equation}
\label{hbound3}
\frac{\partial^{i}}{\partial x^i}\frac{\partial^{j}}{\partial y^j}h(x,y)\ll_{i,j,N} x^{-1-i-j}\left(x^N+\min\{1,(x/|y|)^N\}\right). 
\end{equation}\\

For notational simplicity let us set
$$
r_1=q_1/B,\;\;r_2=q_2/\sqrt{B},\;\;\mathbf{u}=B\mathbf{m}/q_1q_2
$$
so that
$$
I_{q_1,q_2}(\mathbf{m})=\int_{\mathbb R^n}h\left(r_1,Q_1(\mathbf{y})\right)h\left(r_2,r_1^{-1}Q_2(\mathbf{y})\right)U(r_1^{-1}Q_2(\mathbf{y}))
W\left(\mathbf{y}\right)e\left(-\mathbf{u}.\mathbf{y}\right)\mathrm{d}\mathbf{y}.
$$\\

\begin{lemma}
\label{lemma9}
We have
$$
I_{q_1,q_2}(\mathbf{m})\ll_N \left\{\frac{H+r_1^{-1}+(r_1r_2)^{-1}}{|\mathbf{u}|}\right\}^N.
$$\\
\end{lemma}

To prove this lemma one integrates by parts $N$ times, and uses the last bound \eqref{hbound3} for the partial derivatives of $h(x,y)$. Now since $q_1\ll B$ and $q_2\ll B^{1/2}$ (see Lemma~\ref{lemma2}) it follows that the integral in \eqref{integral} is negligibly small if 
$$
|\mathbf{m}|\gg HB^{1/2+\varepsilon}.
$$
So we can cut the sum over $\mathbf{m}$ in \eqref{poisson} at $|\mathbf{m}|\ll HB^{1/2+\varepsilon}$, as the tail makes a negligible contribution. We now proceed towards a detailed analysis of the integral. \\

Let 
$$
f_i(v)=f(r_i,v)=r_i^{k+1}v^j\frac{\partial^{k+j}}{\partial r_i^k\partial v^j}h(r_i,v)
$$
with $k,j=0$ or $1$, and consider
$$
I(\mathbf{u})=\int_{\mathbb R^n}f_1\left(Q_1(\mathbf{y})\right)f_2\left(r_1^{-1}Q_2(\mathbf{y})\right)U\left(r_1^{-1}Q_2(\mathbf{y})\right)W\left(\mathbf{y}\right)e\left(-\mathbf{u}.\mathbf{y}\right)\mathrm{d}\mathbf{y}.
$$
(If one takes $k,j=0$, then one gets $I(\mathbf{u})=r_1r_2I_{q_1,q_2}(\mathbf{m})$.) In the rest of this section we will prove the following lemma.\\

\begin{lemma}
\label{lemma10}
Let $W$ be as in Theorem~\ref{mthm}, $U$ as in \eqref{intro-u} and let $I(\mathbf{u})$ be as defined above. We have
$$
I(\mathbf{u})\ll |\mathbf{u}|^{-n/2}H^nB^\varepsilon.
$$\\
\end{lemma}
 
We choose a suitable weight function $w:\mathbb{R}\rightarrow\mathbb{R}$, with $v^jw^{(j)}(v)\ll_j 1$ such that $w(Q_1(\mathbf{x}))=1$ for $\mathbf{x}\in \mathrm{Supp}(W)$. In particular, we can choose $w$ so that $w(x)=1$ for $x\in [-c',c']$ and $w(x)=0$ for $x\notin [-2c',2c']$ for some $c'>1$ depending only on the support of $W$ and the form $Q_1$. Let
$$
p_1(t)=\int w(v)f_1(v)e(-tv)\mathrm{d}v
$$
and let 
$$
p_2(t)=\int U(v)f_2(v)e(-tv)\mathrm{d}v.
$$
The Fourier transform $p_1(t)$ has been investigated in Section~7 of \cite{H}. We will recall the relevant estimates. Since 
$$
\frac{\partial^N}{\partial v^N}\left(w(v)f_1(v)\right)\ll_N r_1^{-N}\min\{1,(r_1/|v|)^2\},
$$
by repeated integration by parts we get
$$
p_1(t)\ll r_1(r_1|t|)^{-N}
$$
for any $N\geq 0$. In particular $p_1(t)\ll \min\{1,|t|^{-1}\}$. The same bound holds for $p_2$, as we are just replacing $w$ by a similar function $U$.\\

Next using Fourier inversion we get
\begin{align}
\label{exp-int}
I(\mathbf{u})=\iint_{\mathbb{R}^2}p_1(t_1)p_2(t_2)\int_{\mathbb R^n}W(\mathbf{y})e\left(t_1Q_1(\mathbf{y})+t_2r_1^{-1}Q_2(\mathbf{y})-\mathbf{u}.\mathbf{y}\right)\mathrm{d}\mathbf{y}\mathrm{d}t_1\mathrm{d}t_2.
\end{align} 
The inner integral over $\mathbf{y}$ in \eqref{exp-int} is an exponential integral. To estimate it we break the domain of integration into smaller domains, and in each domain either the integral is negligibly small or the trivial bound suffices. To this end let 
$$
w_0(x)=\begin{cases}e^{-\frac{1}{1-x^2}} &\text{if}\;\;|x|\leq 1\\ 0 &\text{otherwise}\end{cases}
$$
and define $c_0=\int w_0(x)\mathrm{d}x$. For any $\delta\in (0,1]$ set 
$$
w_{\delta}(\mathbf{u},\mathbf{v})=c_0^{-n}W(\delta\mathbf{u}+\mathbf{v})\prod_{i=1}^nw_0(u_i),
$$
where $\mathbf{u}=(u_1,\dots,u_n)$. Then
$$
W(\mathbf{y})=\delta^{-n}\int_{\mathbb{R}^n}w_\delta\left(\frac{\mathbf{y}-\mathbf{z}}{\delta},\mathbf{z}\right)\mathrm{d}\mathbf{z}.
$$
The function $w_\delta$ is supported in $[-1,1]^n\times [-c,c]^n$ for some constant $c$, and the function 
$$
G(\mathbf{y})=w_\delta(\delta^{-1}(\mathbf{y}-\mathbf{z}),\mathbf{z})
$$ 
has $\mathrm{Supp}(G)\subset \mathrm{Supp}(W)$ for any given $\mathbf{z}$. Substituting and making a change of variables we arrive at 
\begin{align}
\label{exp-int-bd}
I(\mathbf{u})\leq \int_{[-c,c]^n}\iint_{\mathbb{R}^2}|p_1(t_1)||p_2(t_2)|\Bigl|\int_{\mathbb R^n}w_\delta\left(\mathbf{y},\mathbf{z}\right)e\left(F(\mathbf{y})\right)\mathrm{d}\mathbf{y}\Bigr|\mathrm{d}t_1\mathrm{d}t_2\mathrm{d}\mathbf{z}.
\end{align}
where
$$
F(\mathbf{y})=t_1Q_1(\delta\mathbf{y}+\mathbf{z})+t_2r_1^{-1}Q_2(\delta\mathbf{y}+\mathbf{z})-\mathbf{u}.(\delta\mathbf{y}+\mathbf{z}).
$$\\

To estimate the exponential integral in \eqref{exp-int-bd} we will use the following estimate (see Lemma~10 of \cite{H}). Let $v$ be a compactly supported smooth function on $\mathbb R^n$ with $v^{(\mathbf{j})}\ll_{\mathbf{j}} 1$, and let $f$ be a real valued smooth function on $\mathbb R^n$. Suppose we have positive real numbers $\lambda$ and $\{A_2,A_3,\dots\}$ such that for $\mathbf{x}\in \mathrm{Supp}(v)$ we have $|\nabla f|\geq \lambda$ and
$|v^{(\mathbf{j})}(\mathbf{x})|\leq A_j\lambda$ for $j=j_1+\dots+j_n\geq 2$. Then we have 
\begin{align}
\label{exp-int-est}
\int v(\mathbf{x})e(f(\mathbf{x}))\mathrm{d}\mathbf{x} \ll_{A_i,N} \lambda^{-N}
\end{align}
for any $N\geq 0$. \\

Now we turn our attention at the exponential integral that appears in \eqref{exp-int-bd}. Let $M(t_1,t_2)=t_1M_1+t_2r_1^{-1}M_2$, and let
$$
\tau=\tau(t_1,t_2)=\|M(t_1,t_2)\|_2=\max_{\mathbf{x}\neq \mathbf{0}}\frac{|M(t_1,t_2)\mathbf{x}|}{|\mathbf{x}|}
$$ 
be the Euclidean norm (or the spectral norm) of the matrix. The partial derivatives of $F$ of order $k\geq 2$ are $O(\tau\delta^2)$. If $\nabla F(\mathbf{0})\gg \tau\delta^2B^\varepsilon$ then we have
$$
\nabla F(\mathbf{y})=\nabla F(\mathbf{0})+O(\tau\delta^2)\gg \tau\delta^2B^\varepsilon.
$$  
Suppose $\delta<H^{-1}$, then $w_\delta^{(\mathbf{j})}(\star,\mathbf{z})\ll_\mathbf{j} 1$ for any given $\mathbf{z}$, and by \eqref{exp-int-est} it follows that the inner integral in \eqref{exp-int-bd} is negligibly small whenever 
$$
\nabla F(\mathbf{0})\gg \max\{\tau\delta^2,1\}B^\varepsilon.
$$
Now we consider $|\mathbf{u}|\gg H^2B^\varepsilon$ and pick $\delta=|\mathbf{u}|^{-1/2}$. We say that $(\mathbf{z},t_1,t_2)$ is `good' if 
$$
|\nabla F(\mathbf{0})|=|\mathbf{u}|^{-1/2}|t_1\nabla Q_1(\mathbf{z})+t_2r_1^{-1}\nabla Q_2(\mathbf{z})-\mathbf{u}|\geq c_\varepsilon \max\{\tau/|\mathbf{u}|,1\}B^\varepsilon.
$$
The total contribution of this set is negligibly small. On the other hand if $\tau\leq c|\mathbf{u}|$ for a small constant $c$, then by repeated integration by parts we can show that the integral is small.\\

Our job is now reduced to estimating the size of the `bad' set
$$
\mathfrak{B}(t_1,t_2)=\left\{\mathbf{z}\in [-c,c]^n:|2M(t_1,t_2)\mathbf{z}-\mathbf{u}|< c_\varepsilon B^\varepsilon\frac{\tau}{\sqrt{|\mathbf{u}|}}\right\},
$$
for $\mathbf{u}\neq \mathbf{0}$. If $\tau(t_1,t_2)=0$ (or equivalently $M(t_1,t_2)=0$) then the set is empty and we have $\text{meas}(\mathfrak{B}(t_1,t_2))=0$.
Now consider the case where $\tau> 0$. Suppose $\mathbf{z}$ and $\mathbf{z}+\mathbf{z}'$ are two points in $\mathfrak{B}(t_1,t_2)$. Then 
$$
|M(t_1,t_2)\mathbf{z}'|\leq 10c_\varepsilon B^\varepsilon\frac{\tau}{\sqrt{|\mathbf{u}|}}
$$
and it follows that
$$
\mathrm{meas}(\mathfrak{B}(t_1,t_2))\leq \mathrm{meas}\left\{\mathbf{z}': |M(t_1,t_2)\mathbf{z}'|\leq 10c_\varepsilon B^\varepsilon\frac{\tau}{\sqrt{|\mathbf{u}|}}\right\}.
$$
Since $|M(t_1,t_2)\mathbf{z}'|\leq \tau(t_1,t_2)|\mathbf{z}'|$ the last set is contained in the ball
$$
\left\{\mathbf{z}': |\mathbf{z}'|\leq 10c_\varepsilon B^\varepsilon\frac{1}{\sqrt{|\mathbf{u}|}}\right\}.
$$
Hence we get
$$
\mathrm{meas}(\mathfrak{B}(t_1,t_2))\ll B^\varepsilon|\mathbf{u}|^{-n/2}.
$$ \\

For $\mathbf{z}\in \mathfrak{B}(t_1,t_2)$ we use the trivial bound
\begin{align*}
\int_{\mathbb R^n}w_\delta\left(\mathbf{y},\mathbf{z}\right)e\left(F(\mathbf{y})\right)\mathrm{d}\mathbf{y}\ll 1
\end{align*}
(which follows from the support of $w_\delta(.,\mathbf{z})$) in \eqref{exp-int-bd} to get
\begin{align*}
I(\mathbf{u})\ll \frac{B^\varepsilon}{|\mathbf{u}|^{n/2}}\iint_{\mathbb{R}^2}|p_1(t_1)||p_2(t_2)|\mathrm{d}t_1\mathrm{d}t_2.
\end{align*}
Next using the bounds for the Fourier transforms $p_i(t)$ we get
$$
I(\mathbf{u})\ll \frac{B^\varepsilon}{|\mathbf{u}|^{n/2}}\mathop{\iint}_{\substack{|t_1|\ll r_1^{-1}B^\varepsilon\\|t_2|\ll r_2^{-1}B^\varepsilon}}\min\{1,|t_1|^{-1}\}\min\{1,|t_2|^{-1}\}\mathrm{d}t_1\mathrm{d}t_2 + B^{-2013n}.
$$
We conclude that
$
I(\mathbf{u})\ll |\mathbf{u}|^{-n/2}B^\varepsilon
$
if $|\mathbf{u}|\geq H^2B^\varepsilon$. For smaller values of $|\mathbf{u}|$ we use the trivial bound $I(\mathbf{u})\ll B^\varepsilon$. So it follows that for any $\mathbf{u}$ we have
$$
I(\mathbf{u})\ll |\mathbf{u}|^{-n/2}H^nB^\varepsilon.
$$
This completes the proof of Lemma~\ref{lemma10}.\\

As a consequence of Lemma~\ref{lemma10} we get the following bounds, which will be used in our estimation of the error term in the next section.\\

\begin{lemma}
\label{lemma11}
We have
\begin{align}
\label{int-11}
I_{q_1,q_2}(\mathbf{m})\ll \frac{B^{3/2+\varepsilon}}{q_1q_2}\left(\frac{B|\mathbf{m}|}{q_1q_2}\right)^{-n/2}H^n
\end{align}
and
\begin{align}
\label{int-22}
\frac{\partial}{\partial q_2} I_{q_1,d_2q_2}(\mathbf{m})\ll \frac{B^{3/2+\varepsilon}}{q_1q_2^2d_2}\left(\frac{B|\mathbf{m}|}{q_1q_2d_2}\right)^{-n/2}H^n.
\end{align}\\
\end{lemma}

To prove the first statement we take $k=j=0$ in the definition of $f_i$. For the second statement we first make of change of variables to get
\begin{align*}
I_{q_1,q_2}\left(\mathbf{m}\right)=&q_2^n\int_{\mathbb R^n}h\left(r_1,q_2^2Q_1(\mathbf{y})\right)h\left(r_2,r_1^{-1}q_2^2Q_2(\mathbf{y})\right)\\
&\times U\left(r_1^{-1}q_2^2Q_2(\mathbf{y})\right)W\left(q_2\mathbf{y}\right)e\left(-\frac{B\mathbf{m}.\mathbf{y}}{q_1}\right)\mathrm{d}\mathbf{y}.
\end{align*}
Differentiating with respect to $q_2$, we see that $q_1q_2^2\partial I_{q_1,q_2}(\mathbf{u})/\partial q_2$ can be written as a linear combination of integrals of the type $I(\mathbf{u})$ (with various choices of the pair $(k,j)$). Then we just use the bound from Lemma~\ref{lemma10} for each component. 
\\


\section{The error term}

For notational simplicity we will provide details only for $n$ odd ($n=11$ being the most interesting case). So far we have proved (see Lemma~\ref{lemma2}) that
\begin{align*}
N^\star(B)=B^{n-3}\sum_{q_1\ll B}\;\sum_{q_2\ll B^{1/2}}\;\frac{1}{(q_1q_2)^n}\mathop{\sum\dots\sum}_{|\mathbf{m}|\ll HB^{1/2+\varepsilon}}\;C_{q_1,q_2}(\mathbf{m})\;I_{q_1,q_2}(\mathbf{m}).
\end{align*}
It is expected that the non-zero frequencies contribute to the error term, and accordingly we set
\begin{align}
\label{eb-error}
E(B)=B^{n-3}\sum_{q_1\ll B}\;\sum_{q_2\ll B^{1/2}}\;\frac{1}{(q_1q_2)^n}\mathop{\sum\dots\sum}_{0<|\mathbf{m}|\ll HB^{1/2+\varepsilon}}\;C_{q_1,q_2}(\mathbf{m})\;I_{q_1,q_2}(\mathbf{m}).
\end{align}
In this section we will prove the following estimate.\\

\begin{proposition}
\label{prop1}
We have
\begin{align}
\label{error-bd}
E(B)\ll H^{2n}B^{3n/4-3/2+7/32+\varepsilon}.
\end{align}\\
\end{proposition}

Taking absolute values and using mulitplicativity of the character sum we get
\begin{align*}
E(B) \leq B^{n-3}\sum_{q_1\ll B}\;\frac{1}{q_1^n}&\mathop{\sum\dots\sum}_{0<|\mathbf{m}|\ll HB^{1/2+\varepsilon}}\;\sum_{d_2|(q_1N(\mathbf{m}))^\infty}\frac{|C_{q_1,d_2}(\mathbf{m})|}{d_2^n}\\
&\times \Bigl|\sum_{\substack{q_2\ll B^{1/2}/d_2\\(q_2,q_1N(\mathbf{m}))=1}} C_{1,q_2}(\mathbf{m})\;q_2^{-n}I_{q_1,q_2d_2}(\mathbf{m})\Bigr|
\end{align*}
where $N(\mathbf{m})$ is as defined in \eqref{nm}. To the inner sum over $q_2$ we apply partial summation to get
$$
\sum_{\substack{q_2\ll B^{1/2}/d_2\\(q_2,q_1N(\mathbf{m}))=1}}\: C_{1,q_2}(\mathbf{m})\;q_2^{-n}I_{q_1,q_2d_2}(\mathbf{m})=-\int_1^{B^{1/2}/d_2} \Bigl[\sum_{\substack{q_2\leq x\\(q_2,q_1N(\mathbf{m}))=1}} C_{1,q_2}(\mathbf{m})\Bigr]\frac{\partial}{\partial x}\frac{I_{q_1,d_2 x}(\mathbf{m})}{x^n}\mathrm{d}x.
$$
To the sum over $q_2$ we use the bounds from \eqref{char-sum-bd-1} and \eqref{char-sum-bd-0}, and to the derivative we apply the bounds \eqref{int-11} and \eqref{int-22}, to get
$$
\sum_{\substack{q_2\ll B^{1/2}/d_2\\(q_2,q_1N(\mathbf{m}))=1}} C_{1,q_2}(\mathbf{m})\;q_2^{-n}I_{q_1,q_2d_2}(\mathbf{m})\ll \frac{H^nB^{3/2+\varepsilon}|
\mathbf{m}|^{\theta_1(\mathbf{m})}}{q_1d_2}\left(\frac{q_1d_2}{B|\mathbf{m}|}\right)^{n/2}\left(\frac{B}{d_2}\right)^{\theta_2(\mathbf{m})}
$$
where 
$$
\theta_1(\mathbf{m})=\begin{cases}7/16 &\text{if}\;\;(-1)^{(n-1)/2}Q_2^\star(\mathbf{m})\neq \Box\\
0&\text{otherwise},
\end{cases}
$$ 
(recall that we are assuming that $n$ is odd), and 
$$
\theta_2(\mathbf{m})=\begin{cases} 0 & \text{if} \;\;(-1)^{(n-1)/2}Q_2^\star(\mathbf{m})\neq \Box\\
1/4&\text{otherwise}.\end{cases}
$$ 
Moreover from Lemma~\ref{lemma4} it follows that
$$
\sum_{\substack{d_2|N(\mathbf{m})^\infty\\d_2\ll B^{1/2}}}\frac{|C_{1,d_2}(\mathbf{m})|}{d_2^{n/2+1}}\ll B^\varepsilon
$$
as $0\neq N(\mathbf{m})\ll H^2B^{1+\varepsilon}$.\\

Hence we get
\begin{align*}
E(B)\ll H^nB^{n/2-3/2+\varepsilon}&\mathop{\sum\dots\sum}_{0<|\mathbf{m}|\ll HB^{1/2+\varepsilon}}|\mathbf{m}|^{\theta_1(\mathbf{m})-n/2}\mathop{\sum\sum}_{\substack{q_1\ll B\\d_2\ll B^{1/2}\\ d_2|q_1^\infty}}\frac{|C_{q_1,d_2}(\mathbf{m})|}{(q_1d_2)^{n/2+1}}\left(\frac{B}{d_2}\right)^{\theta_2(\mathbf{m})}.
\end{align*}
The last sum over $q_1$ and $d_2$, is dominated by
\begin{align*}
\mathop{\sum\sum}_{\substack{d_1\ll B\\ d_2\ll B^{1/2}\\d_2|d_1^\infty\\d_1|(d_2g(\mathbf{m}))^\infty}}\;\frac{|C_{d_1,d_2}(\mathbf{m})|}{(d_1d_2)^{n/2+1}}\left(\frac{B}{d_2}\right)^{\theta_2(\mathbf{m})}\sum_{\substack{q_1\ll B/d_1\\(q_1,g(\mathbf{m}))=1}}\;\frac{|C_{q_1,1}(\mathbf{m})|}{q_1^{n/2+1}},
\end{align*}
where we set $g(\mathbf{m})=\Delta_VG(\mathbf{m})$ if $G(\mathbf{m})\neq 0$ and $g(\mathbf{m})=\Delta_V\text{gcd}(\mathbf{m})$ if $G(\mathbf{m})= 0$. Applying partial summation together with Lemma~\ref{lemma7}, we arrive at
\begin{align*}
B^\varepsilon \mathop{\sum\sum}_{\substack{d_1\ll B\\ d_2\ll B^{1/2}\\d_2|d_1^\infty\\d_1|(d_2g(\mathbf{m}))^\infty}}\;\frac{|C_{d_1,d_2}(\mathbf{m})|}{(d_1d_2)^{n/2+1}}\left(\frac{B}{d_2}\right)^{\theta_2(\mathbf{m})}\left(\frac{B}{d_1} \right)^{\delta(G(\mathbf{m}))/2}.
\end{align*}
Using Lemma~\ref{lemma6} and Lemma~\ref{lemma8}, we bound the above sum by 
\begin{align}
\label{whatnext}
B^{\varepsilon} \mathop{\sum\sum}_{\substack{d_1\ll B\\ d_2\ll B^{1/2}\\d_2|d_1^\infty\\d_1|d_2^\infty}}\;\frac{\left(\text{rad}(d_2),Q_2^\star(\mathbf{m})\right)}{\text{rad}(d_2)}\left(\frac{B}{d_2}\right)^{\theta_2(\mathbf{m})}\left(\frac{B}{d_1} \right)^{\delta(G(\mathbf{m}))/2}.
\end{align}
We need to consider four distinct situations. Suppose $G(\mathbf{m})=Q_2^\star(\mathbf{m})=0$ then the sum boils down to
\begin{align*}
B^{3/4+\varepsilon} \mathop{\sum\sum}_{\substack{d_1\ll B\\ d_2\ll B^{1/2}\\d_2|d_1^\infty\\d_1|d_2^\infty}}\;\frac{1}{\sqrt{d_1}d_2^{1/4}}\ll B^{1+\varepsilon}.
\end{align*}
If $G(\mathbf{m})=0$ and $Q_2^\star(\mathbf{m})\neq 0$ then the sum in \eqref{whatnext} is bounded by
\begin{align*}
B^{\theta_2(\mathbf{m})+1/2+\varepsilon} \mathop{\sum\sum}_{\substack{d_1\ll B\\ d_2\ll B^{1/2}\\d_2|d_1^\infty\\d_1|d_2^\infty}}\;\frac{\left(\text{rad}(d_2),Q_2^\star(\mathbf{m})\right)}{\text{rad}(d_2)}\ll B^{\theta_2(\mathbf{m})+1/2+\varepsilon}\ll B^{1+\varepsilon}.
\end{align*}
On the other hand if $G(\mathbf{m})\neq 0$ and $Q_2^\star(\mathbf{m})= 0$ then we get
\begin{align*}
B^{1/4+\varepsilon} \mathop{\sum\sum}_{\substack{d_1\ll B\\ d_2\ll B^{1/2}\\d_2|d_1^\infty\\d_1|d_2^\infty}}\;\frac{1}{d_2^{1/4}}\ll B^{1/4+3/8+\varepsilon}\ll B^{1+\varepsilon}.
\end{align*}
Finally if $G(\mathbf{m})Q_2^\star(\mathbf{m})\neq 0$ then the sum in \eqref{whatnext} is bounded by
\begin{align*}
B^{\theta_2(\mathbf{m})+\varepsilon} \mathop{\sum\sum}_{\substack{d_1\ll B\\ d_2\ll B^{1/2}\\d_2|d_1^\infty\\d_1|d_2^\infty}}\;\frac{\left(\text{rad}(d_2),Q_2^\star(\mathbf{m})\right)}{\text{rad}(d_2)}\ll B^{\theta_2(\mathbf{m})+\varepsilon}.
\end{align*}
It follows that
\begin{align*}
E(B)\ll H^nB^{n/2-3/2+\varepsilon}&\mathop{\sum\dots\sum}_{0<|\mathbf{m}|\ll HB^{1/2+\varepsilon}}|\mathbf{m}|^{\theta_1(\mathbf{m})-n/2}B^{\psi(\mathbf{m})},
\end{align*}
where $\psi(\mathbf{m})=1$ if $G(\mathbf{m})Q_2^\star(\mathbf{m})=0$ and $\psi(\mathbf{m})=\theta_2(\mathbf{m})$ otherwise.\\

First consider the case where $Q_2^\star(\mathbf{m})G(\mathbf{m})=0$. The contribution of this part to $E(B)$ is bounded by
\begin{align*}
H^{n+7/16}B^{n/2-1/2+7/32+\varepsilon}\mathop{\sum\dots\sum}_{\substack{0<|\mathbf{m}|\ll HB^{1/2+\varepsilon}\\Q_2^\star(\mathbf{m})G(\mathbf{m})=0}}|\mathbf{m}|^{-n/2},
\end{align*}
which is dominated by
\begin{align*}
H^{n+7/16}B^{n/2-1/2+7/32+\varepsilon}\mathop{\max}_{1/2<M\ll HB^{1/2+\varepsilon}}M^{-n/2}\mathop{\sum\dots\sum}_{\substack{0<|\mathbf{m}|\ll M\\Q_2^\star(\mathbf{m})G(\mathbf{m})=0}}1\ll H^{2n}B^{3n/4-3/2+7/32+\varepsilon}.
\end{align*}
(To estimate the last sum we use (7.1) from \cite{BM1}.)
We are little wasteful in the matter of the power of $H$, but this is of little consequence. Next consider the case where $Q_2^\star(\mathbf{m})G(\mathbf{m})\neq 0$, but $(-1)^{(n-1)/2}Q_2^\star(\mathbf{m})=\Box$. The contribution of this part to $E(B)$ is dominated by
\begin{align*}
H^{n}B^{n/2-5/4+\varepsilon}\mathop{\sum\dots\sum}_{\substack{0<|\mathbf{m}|\ll HB^{1/2+\varepsilon}\\(-1)^{(n-1)/2}Q_2^\star(\mathbf{m})=\Box}}|\mathbf{m}|^{-n/2},
\end{align*}
which is dominated by (see (7.2) of \cite{BM1})
\begin{align*}
H^{2n}B^{3n/4-7/4+\varepsilon}.
\end{align*}
Finally the contribution of the generic case, where $G(\mathbf{m})\neq 0$ and $(-1)^{(n-1)/2}Q_2^\star(\mathbf{m})\neq \Box$, to $E(B)$ is dominated by
\begin{align*}
H^{n+7/16}B^{n/2-1/2+7/32+\varepsilon}\mathop{\sum\dots\sum}_{\substack{0<|\mathbf{m}|\ll HB^{1/2+\varepsilon}}}|\mathbf{m}|^{-n/2}\ll H^{2n}B^{3n/4-3/2+7/32+\varepsilon}.
\end{align*}
This completes the proof of the proposition.\\


\section{The main term}
\label{maintermsection}

In this section we will compute the contribution of the zero frequency which is given by
\begin{align*}
M(B)=B^{n-3}\sum_{q_1\ll B}\;\sum_{q_2\ll B^{1/2}}\;\frac{1}{(q_1q_2)^n}\;C_{q_1,q_2}(\mathbf{0})\;I_{q_1,q_2}(\mathbf{0}).
\end{align*}
First we note that from Lemmas~\ref{lemma4}, \ref{lemma6} and \ref{lemma8}, it follows that 
$
C_{q_1,q_2}(\mathbf{0})\ll (q_1q_2)^{n/2+1+\varepsilon}.
$
Also a trivial estimation yields $I_{q_1,q_2}(\mathbf{0})\ll B^{3/2}(q_1q_2)^{-1}$. Consequently we get
\begin{align*}
B^{n-3}\sum_{q_1\ll B}\;\sum_{q_2\geq B^{1/2-1/4n}}\;\frac{1}{(q_1q_2)^n}\;C_{q_1,q_2}(\mathbf{0})\;I_{q_1,q_2}(\mathbf{0})\ll B^{3n/4-3/2+1/8+\varepsilon}.
\end{align*}
For $q_2<B^{1/2-1/4n}$ it follows from \eqref{hbound3} that 
$$
h\left(\frac{q_2}{\sqrt{B}},\frac{BQ_2(\mathbf{y})}{q_1}\right)
$$
is negligibly small if $B|Q_2(\mathbf{y})|q_1^{-1}>1/2$. Consequently we get
\begin{align*}
I_{q_1,q_2}(\mathbf{0})=\tilde{I}_{q_1,q_2}(\mathbf{0})+O(B^{-2013})
\end{align*}
if $q_2<B^{1/2-1/4n}$, where
\begin{align}
\label{integral-tilde}
\tilde{I}_{q_1,q_2}(\mathbf{0})=&\int_{\mathbb R^n}h\left(\frac{q_1}{B},Q_1(\mathbf{y})\right)h\left(\frac{q_2}{\sqrt{B}},\frac{BQ_2(\mathbf{y})}{q_1}\right) W\left(\mathbf{y}\right)\mathrm{d}\mathbf{y}.
\end{align}
Hence we get
\begin{align*}
M(B)=B^{n-3}\sum_{q_1\ll B}\;\sum_{q_2< B^{1/2-1/4n}}\;\frac{1}{(q_1q_2)^n}\;C_{q_1,q_2}(\mathbf{0})\;\tilde{I}_{q_1,q_2}(\mathbf{0})+O(B^{3n/4-3/2+1/8+\varepsilon}).
\end{align*}
Using the above estimate for the character sum, and the bound $\tilde{I}_{q_1,q_2}(\mathbf{0})\ll B^{3/2}(q_1q_2)^{-1}$, we can now complete the sums over $q_1$ and $q_2$ without worsening the error term. \\

\begin{lemma}
\label{lemma12}
We have
$$
M(B)=M^\star(B)+O(B^{3n/4-3/2+1/8+\varepsilon})
$$
where
\begin{align*}
M^\star(B)=B^{n-3}\sum_{q_1=1}^\infty\;\sum_{q_2=1}^{\infty}\;\frac{1}{(q_1q_2)^n}\;C_{q_1,q_2}(\mathbf{0})\;\tilde{I}_{q_1,q_2}(\mathbf{0}).
\end{align*}
with $\tilde{I}_{q_1,q_2}$ as given in \eqref{integral-tilde}.\\
\end{lemma}

Next consider the double Dirichlet series
\begin{align}
\label{dirserd}
D(s_1,s_2)=\sum_{q_1=1}^\infty\;\sum_{q_2=1}^{\infty}\;\frac{C_{q_1,q_2}(\mathbf{0})}{q_1^{s_1}q_2^{s_2}}.
\end{align}
This is given by an Euler product where the factor corresponding to the prime $p$ is given by
\begin{align*}
D_p(s_1,s_2)&=\sum_{r=0}^\infty\;\sum_{\ell=0}^{\infty}\;\frac{C_{p^r,p^\ell}(\mathbf{0})}{p^{rs_1+\ell s_2}}.
\end{align*}
The following lemma follows from the bound $C_{q_1,q_2}(\mathbf{0})\ll (q_1q_2)^{n/2+1+\varepsilon}$ we noted above.\\

\begin{lemma}
\label{lemma13}
The Dirichlet series $D(s_1,s_2)$, as defined in \eqref{dirserd}, converges absolutely in the domain 
$$
\sigma_1>n/2+2+\varepsilon,\;\;\;\text{and}\;\;\; \sigma_2>n/2+2+\varepsilon.
$$
Moreover in this region we have $D(s_1,s_2)\ll 1$. \\
\end{lemma}

Theorem~\ref{mthm} will follow from Proposition~\ref{prop1} and the following asymptotic, which will be proved in the rest of this section.\\

\begin{proposition}
\label{prop2}
Suppose $n>6$ and $H\geq 1$, then we have
\begin{align}
\label{maintermeqn}
M(B)=B^{n-4}D(n-1,n)J_0(W)+O(H^2B^{n-5+\varepsilon}+HB^{3n/4-3/2+1/8+\varepsilon}),
\end{align}
where $J_0(W)$ is as defined in \eqref{j0w}.\\
\end{proposition}

Let $W:\mathbb{R}^n\rightarrow \mathbb{R}$ be as in Theorem~\ref{mthm}. We want to analyse the integral
\begin{align}
\label{integralv}
W(s_1,s_2)=\int_{\mathbb R^n}|Q_1(\mathbf{y})|^{s_1-1}|Q_2(\mathbf{y})|^{s_2-1} W\left(\mathbf{y}\right)\mathrm{d}\mathbf{y},
\end{align}
which we consider as a function of two complex variables. A priori, the integral is defined in the tube domain $\sigma_1> 1$ and $\sigma_2> 1$, where it is a holomorphic function in two variables. Now we will show that the function has a meromorphic extension to $\mathbb{C}^2$. 
Recall that $V^\star$ denotes the singular locus, which is the union of the eigenspaces of $M_2^{-1}M_1$. Since we are assuming that $V^\star\cap \text{Supp}(W)=\emptyset$, for every $\mathbf{z}\in \text{Supp}(W)$ the matrix 
$$
(M_1\mathbf{z}\: M_2\mathbf{z})
$$ 
is of rank two. For any $1\leq i<j\leq n$ we define $\mathbf{y}_{i,j}$ to be the vector which is obtained from $\mathbf{y}$ by deleting the $i$-th and the $j$-th entry. Let
\begin{align}
\label{phiij}
\phi_{i,j}:\mathbb{R}^n\mapsto\mathbb{R}^n,\;\;\;\phi_{i,j}(\mathbf{y})=(Q_1(\mathbf{y}), Q_2(\mathbf{y}),\mathbf{y}_{i,j}),
\end{align}
and let $J_{\phi_{i,j}}$ be the associated Jacobian. Consider the map
$$
\mathbf{y}\mapsto (J_{\phi_{i,j}}(\mathbf{y}))_{1\leq i<j\leq n}.
$$
This map takes the set $\mathrm{Supp}(W)$ to a compact set which does not contain $\mathbf{0}$. Consequently there is a $\theta>0$, depending only on $W$, such  that 
$$
\max_{i,j} |J_{\phi_{i,j}}(\mathbf{y})|\geq \theta
$$
for all $\mathbf{y}\in\mathrm{Supp}(W)$. So it follows that we have a smooth finite partition 
$$
W(\mathbf{y})=\mathop{\sum\sum}_{1\leq i<j\leq n} \sum_{\xi\in \Xi(i,j)} W_{i,j}^\xi(\mathbf{y})
$$
such that $\mathrm{Supp}(W_{i,j}^\xi)$ is connected and for any $\mathbf{y}$ in this set we have $|J_{\phi_{i,j}}(\mathbf{y})|\geq \theta/2$. In particular in $\mathrm{Supp}(W_{i,j}^\xi)$ the map $\phi_{i,j}$ is one-to-one. Let $\psi_{i,j}^\xi$ be its inverse and let $J_{i,j}^\xi$ be the Jacobian of this map. By the smooth version of the inverse function theorem we have that $\psi_{i,j}^\xi$ is smooth.  \\


We have
$$
W(s_1,s_2)=\mathop{\sum\sum}_{1\leq i<j\leq n} \sum_{\xi\in \Xi(i,j)} W_{i,j}^\xi(s_1,s_2)
$$
where
\begin{align*}
W_{i,j}^\xi(s_1,s_2)&=\int_{\mathbb R^n}|Q_1(\mathbf{y})|^{s_1-1}|Q_2(\mathbf{y})|^{s_2-1} W_{i,j}^\xi\left(\mathbf{y}\right)\mathrm{d}\mathbf{y}\\
&=\int_{\mathbb R^2}|u_1|^{s_1-1}|u_2|^{s_2-1}\mathcal{W}_{i,j}^\xi(u_1,u_2)\mathrm{d}\mathbf{u} 
\end{align*}
and
\begin{align}
\label{wijxi}
\mathcal{W}_{i,j}^\xi(u_1,u_2)&=\int_{\mathbb{R}^{n-2}}W_{i,j}^\xi\left(\psi_{i,j}^\xi(u_1,u_2,\mathbf{y}_{i,j})\right)
|J_{i,j}^\xi\left(u_1,u_2,\mathbf{y}_{i,j}\right)|\mathrm{d}\mathbf{y}_{i,j}\\
\nonumber &=\int_{\mathbb{R}^{n-2}}W_{i,j}^\xi\left(\psi_{i,j}^\xi(u_1,u_2,\mathbf{y}_{i,j})\right)
\left|J_{\phi_{i,j}}\left(\psi_{i,j}^\xi(u_1,u_2,\mathbf{y}_{i,j})\right)\right|^{-1}\mathrm{d}\mathbf{y}_{i,j}.
\end{align}
Since the Jacobian $J_{\phi_{i,j}}$ is bounded away from zero in the support of $W_{i,j}^\xi$ and $\psi_{i,j}^\xi$ is smooth, it follows that $\mathcal{W}_{i,j}^\xi(u_1,u_2)$ is smooth. Furthermore
$$
\frac{\partial^{m_1+m_2}}{\partial u_1^{m_1}\partial u_2^{m_2}}\mathcal{W}_{i,j}^\xi(u_1,u_2)\ll H^{m_1+m_2}.
$$
Using integration by parts we can now analytically continue $W(s_1,s_2)$ to whole of $\mathbb{C}^2$.\\

Next we want to compute the double residue
$$
\mathop{\text{Res}}_{\substack{s_1=0\\s_2=0}}W(s_1,s_2).
$$
There is a box $[-\theta,\theta]^2$ where we have the Taylor expansion 
$$
\mathcal{W}_{i,j}^\xi(u_1,u_2)=\mathcal{W}_{i,j}^\xi(0,0)+u_1\frac{\partial}{\partial u_1}\mathcal{W}_{i,j}^\xi(0,0)+u_2\frac{\partial}{\partial u_2}\mathcal{W}_{i,j}^\xi(0,0)+\dots.
$$
Now splitting the integral as
\begin{align*}
W_{i,j}^\xi(s_1,s_2)&=\int_{[-\theta,\theta]^2}+\int_{\mathbb{R}^2-[-\theta,\theta]^2}\;|u_1|^{s_1-1}|u_2|^{s_2-1}\mathcal{W}_{i,j}^\xi(u_1,u_2)\mathrm{d}\mathbf{u},
\end{align*}
we see that the double residue of the second integral vanishes, and using the Taylor expansion we find that the double residue of the first integral is $4\mathcal{W}_{i,j}^\xi(0,0)$. Hence
\begin{align}
\label{wijxi-res}
\mathop{\text{Res}}_{\substack{s_1=0\\s_2=0}}W^\xi_{i,j}(s_1,s_2)=4\mathcal{W}_{i,j}^\xi(0,0).
\end{align}
Also we have
\begin{align*}
\mathcal{W}_{i,j}^\xi(0,0)&=\mathop{\lim\lim}_{\varepsilon_1,\varepsilon_2\rightarrow 0}\frac{1}{4\varepsilon_1\varepsilon_2}\int_{\substack{|u_1|<\varepsilon_1\\|u_2|<\varepsilon_2}}\mathcal{W}_{i,j}^\xi(u_1,u_2)\mathrm{d}\mathbf{u} \\
&=\mathop{\lim\lim}_{\varepsilon_1,\varepsilon_2\rightarrow 0}\frac{1}{4\varepsilon_1\varepsilon_2}\int_{\substack{|Q_1(\mathbf{y})|<\varepsilon_1\\|Q_2(\mathbf{y})|<\varepsilon_2}}W_{i,j}^\xi(\mathbf{y})\mathrm{d}\mathbf{y}.
\end{align*}
Summing over $\xi$ and $i$, $j$, we get
\begin{align*}
\mathop{\text{Res}}_{\substack{s_1=0\\s_2=0}}W(s_1,s_2)&=4\mathop{\sum\sum}_{1\leq i<j\leq n} \sum_{\xi\in \Xi(i,j)}\mathcal{W}_{i,j}^\xi(0,0) \\
&=4\mathop{\lim\lim}_{\varepsilon_1,\varepsilon_2\rightarrow 0}\frac{1}{4\varepsilon_1\varepsilon_2}\int_{\substack{|Q_1(\mathbf{y})|<\varepsilon_1\\|Q_2(\mathbf{y})|<\varepsilon_2}}W(\mathbf{y})\mathrm{d}\mathbf{y}.
\end{align*}
One may compare this with Heath-Brown's computation of the singular integral in \cite{H}.\\

\begin{lemma}
\label{lemma14}
The function $W(s_1,s_2)$ has a meromorphic continuation to all of $\mathbb{C}^2$ with possible polar divisors at $s_1, s_2=0,-1,\dots$. Moreover we have
$$
W(s_1,s_2)\ll \frac{H^{m_1+m_2}}{\|s_1\|^{m_1}\|s_2\|^{m_2}},
$$
where $\|s_i\|=\min\{|s_i-n|:n=0,-1,\dots\}$ and $m_i\geq \max\{0,1-\sigma_i\}$. Also
\begin{align}
\label{j0w}
J_0(W)=\mathop{\mathrm{Res}}_{\substack{s_1=0\\s_2=0}}\frac{W(s_1,s_2)}{4}=\mathop{\lim\lim}_{\varepsilon_1,\varepsilon_2\rightarrow 0}\frac{1}{4\varepsilon_1\varepsilon_2}\int_{\substack{|Q_1(\mathbf{y})|<\varepsilon_1\\|Q_2(\mathbf{y})|<\varepsilon_2}}W(\mathbf{y})\mathrm{d}\mathbf{y}.
\end{align}\\
\end{lemma}

We will now analyse the integral $\tilde I_{q_1,q_2}(\mathbf{0})$ using Mellin transform. We start by looking at the Mellin transform of $h(x,y)$ for any given $y\neq 0$. Using the definition \eqref{hxy-form} we get 
\begin{align*}
\int_0^\infty h(x,y)x^{s-1}\mathrm{d}x=\zeta(s)\left[\tilde{w}(s-1)-|y|^{s-1}\tilde{w}(1-s)\right]
\end{align*}
for $\sigma>1$, where $\tilde{w}$ is the Mellin transform of $w$. The right hand side, which we denote by $H(s,y)$, extends to an entire function and is of rapid decay in any vertical strip. \\

Applying Mellin inversion we get
\begin{align*}
\tilde I_{q_1,q_2}(\mathbf{0})=&\frac{1}{(2\pi i)^2}\mathop{\iint}_{(\sigma_i)}\left(\frac{B}{q_1}\right)^{s_1}\left(\frac{\sqrt{B}}{q_2}\right)^{s_2}\\
&\times \int_{\mathbb R^n} H\left(s_1,Q_1(\mathbf{y})\right)H\left(s_2,\frac{BQ_2(\mathbf{y})}{q_1}\right)W\left(\mathbf{y}\right)\mathrm{d}\mathbf{y}\;\mathrm{d}\mathbf{s}
\end{align*}
where for the convergence of the last integral we require to take $\sigma_1, \sigma_2\geq 1$. (Note that the set $Q_i(\mathbf{y})=0$ has measure zero.) Using the above explicit expression for the Mellin transform we get
\begin{align*}
\tilde I_{q_1,q_2}(\mathbf{0})=\mathop{\sum\sum}_{\delta_1,\delta_2=0,1}\;(-1)^{\delta_1+\delta_2}\:I_{q_1,q_2}^{\delta_1,\delta_2},
\end{align*}
where
\begin{align*}
I_{q_1,q_2}^{\delta_1,\delta_2}=&\frac{1}{(2\pi i)^2}\mathop{\iint}_{(\sigma_i)}\left(\frac{B}{q_1}\right)^{s_1+\delta_2(s_2-1)}\left(\frac{\sqrt{B}}{q_2}\right)^{s_2}\zeta(s_1)\zeta(s_2)\\
&\times \tilde{w}\left((-1)^{\delta_1}(s_1-1)\right)\tilde{w}\left((-1)^{\delta_2}(s_2-1)\right)W\left(s_1^{\delta_1},s_2^{\delta_2}\right)\mathrm{d}\mathbf{s}.
\end{align*}\\

Accordingly we get
\begin{align*}
M^\star(B)=\mathop{\sum\sum}_{\delta_1,\delta_2=0,1}\;(-1)^{\delta_1+\delta_2}\:M^{\delta_1,\delta_2}(B),
\end{align*}
where
\begin{align*}
M^{\delta_1,\delta_2}(B)=&\frac{B^{n-3}}{(2\pi i)^2}\mathop{\iint}_{(\sigma_i)}\:B^{s_1+s_2/2+\delta_2(s_2-1)}\zeta(s_1)\zeta(s_2)D(n+s_1+\delta_2(s_2-1),n+s_2)\\
&\times \tilde{w}\left((-1)^{\delta_1}(s_1-1)\right)\tilde{w}\left((-1)^{\delta_2}(s_2-1)\right)W\left(s_1^{\delta_1},s_2^{\delta_2}\right)\mathrm{d}\mathbf{s}.
\end{align*}
Here $\sigma_1,\sigma_2=1+\varepsilon$. Now we want to move the contours to the left. The residues at the poles of the zeta functions cancel pairwise. So we can first move to $\sigma_1,\sigma_2=\varepsilon$. For the terms with $\delta_2=0$, we can shift the contour $(\sigma_2)$ further to the left, upto $\sigma_2=-n/2+2+\varepsilon$. Then using the bounds from Lemmas~\ref{lemma13} and \ref{lemma14} we get
$$
M^{\delta_1,0}(B)\ll HB^{3n/4-2+\varepsilon}.
$$
For the term with $(\delta_1,\delta_2)=(0,1)$, we can shift the contour $(\sigma_1)$ to $\sigma_1=-n/2+3+\varepsilon$. Using the bounds from Lemmas~\ref{lemma13} and \ref{lemma14} we get
$$
M^{0,1}(B)\ll HB^{n/2-1+\varepsilon}.
$$
This is dominated by the previous bound as $n>4$. We conclude that
$$
M^\star(B)=M^{1,1}(B)+O(HB^{3n/4-2+\varepsilon}).
$$\\

We move $\sigma_1$ to $-1+\varepsilon$, collecting the residue at the polar divisor $s_1=0$ and using the bounds from Lemmas~\ref{lemma13} and \ref{lemma14} we obtain
\begin{align*}
M^{1,1}(B)=\frac{-1}{4\pi i}\mathop{\int}_{(\sigma_2)}&B^{n+3s_2/2-4}\zeta(s_2)D(n+s_2-1,n+s_2)\\
&\times\tilde{w}(1-s_2)\mathop{\text{Res}}_{s_1=0} W(s_1,s_2)\mathrm{d}\mathbf{s}+O(H^2B^{n-5+\varepsilon}),
\end{align*}
as $\zeta(0)=-1/2$ and $\tilde{w}(1)=1$. Next we shift $\sigma_2$ to $-1+\varepsilon$, passing through a pole at $s_2=0$, and obtain
\begin{align*}
M^{1,1}(B)=\frac{1}{4}B^{n-4}D(n-1,n)\mathop{\text{Res}}_{\substack{s_1=0\\s_2=0}}W(s_1,s_2)
+O(H^2B^{n-5+\varepsilon}).
\end{align*}
This completes the proof of Proposition~\ref{prop2}.\\


\section{Sharp cuts}
\label{sharp-cuts}

In this section we will prove Theorem~\ref{thm2}. Recall that $\mathcal P=\prod_{i=1}^n[c_i,d_i]$ is a box in $\mathbb R^n$, which does not intersect the singular locus. We can slightly thicken the box $\mathcal{P}$, say by a parameter $\delta=1/H$, to get $\mathcal{P}'$ so that the condition on the support still holds. Now there is a non-negative smooth function $W$ supported in $\mathcal{P}'$, such that $W(\mathbf{y})=1$ for $\mathbf{y}\in\mathcal{P}$, and such that $W^{(j)}\ll H^j$. We can also find a smooth function $V$ supported in $\mathcal{P}$ such that $V(\mathbf{x})=1$ for any $\mathbf{x}\in \mathcal{P}$ which is a distance of $1/H$ away from the set $\mathcal{P}'-\mathcal{P}$. We may also have such a $V$ with $0\leq V\leq 1$ and $V^{(j)}\ll H^j$. \\

Let $U=W-V$, which is a non-negative smooth function with
$$\mathrm{Supp}(U)\subset \prod_{1\leq i\leq n}\left\{[c_i-H^{-1},c_i+H^{-1}]\cup [d_i-H^{-1},d_i+H^{-1}]\right\}.$$ 
Then we have
$$
\mathop{\sum\dots\sum}_{\substack{\mathbf{m}\in B\mathcal{P}\\Q_1(\mathbf{m})=Q_2(\mathbf{m})=0}}1=
\mathop{\sum\dots\sum}_{Q_1(\mathbf{m})=Q_2(\mathbf{m})=0}W\left(\frac{\mathbf{m}}{B}\right)+O\left(
\mathop{\sum\dots\sum}_{Q_1(\mathbf{m})=Q_2(\mathbf{m})=0}U\left(\frac{\mathbf{m}}{B}\right)\right).
$$
Applying Theorem~\ref{mthm} to both the main term and the error term we get
$$
\mathfrak{S}J_0(W)B^{n-4}+
O\left(J_0(U)B^{n-4}+H^2B^{n-5+\varepsilon}+H^{2n}B^{3n/4-41/32+\varepsilon}\right).
$$
We will obtain an upper bound for $J_0(U)$. We return to the expression \eqref{wijxi}, where $U$ now takes the place of $W$. Since 
$$\mathrm{Supp}(U_{i,j}^\xi)\subset \mathrm{Supp}(U),
$$ from \eqref{wijxi} we deduce that
$$
\mathcal{U}_{1,2}^\xi(0,0)\ll \mathop{\int\dots\int}_{\prod_{3\leq i\leq n}\left\{[c_i-H^{-1},c_i+H^{-1}]\cup [d_i-H^{-1},d_i+H^{-1}]\right\}} \mathbf{y}_{1,2}\ll H^{-n+2}.
$$
The same bound holds for all $\mathcal{U}_{i,j}^\xi(0,0)$, and consequently we have
$$
J_0(U)\ll H^{-n+2}.
$$
We conclude that
$$
\mathop{\sum\dots\sum}_{\substack{\mathbf{m}\in B\mathcal{P}\\Q_1(\mathbf{m})=Q_2(\mathbf{m})=0}}1=
\mathfrak{S}J_0(W)B^{n-4}+
O\left(H^{-n+2}B^{n-4}+H^2B^{n-5+\varepsilon}+H^{2n}B^{3n/4-41/32+\varepsilon}\right).
$$
The error term is satisfactory, of the size $O(B^{n-4-\delta})$ for some positive $\delta$, if we pick $H=B^{1/100n}$ and $n\geq 11$.\\

In the main term the singular integral $J_0(W)$ depends on $W$. But from \eqref{j0w} we get
\begin{align*}
J_0(W)&=\mathop{\lim\lim}_{\varepsilon_1,\varepsilon_2\rightarrow 0}\frac{1}{4\varepsilon_1\varepsilon_2}\int_{\substack{|Q_1(\mathbf{y})|<\varepsilon_1\\|Q_2(\mathbf{y})|<\varepsilon_2}}W(\mathbf{y})\mathrm{d}\mathbf{y}\\
&=\mathop{\lim\lim}_{\varepsilon_1,\varepsilon_2\rightarrow 0}\frac{1}{4\varepsilon_1\varepsilon_2}\int_{\substack{|Q_1(\mathbf{y})|<\varepsilon_1\\|Q_2(\mathbf{y})|<\varepsilon_2\\\mathbf{y}\in\mathcal{P}}}\mathrm{d}\mathbf{y}+\mathop{\lim\lim}_{\varepsilon_1,\varepsilon_2\rightarrow 0}\frac{1}{4\varepsilon_1\varepsilon_2}\int_{\substack{|Q_1(\mathbf{y})|<\varepsilon_1\\|Q_2(\mathbf{y})|<\varepsilon_2\\\mathbf{y}\notin\mathcal{P}}}\mathrm{d}\mathbf{y}.
\end{align*} 
The first integral is just $J_0(\mathcal{P})$, the singular integral for the indicator function of the box $\mathcal{P}$, and the second integral is dominated by 
$$
\mathop{\lim\lim}_{\varepsilon_1,\varepsilon_2\rightarrow 0}\frac{1}{4\varepsilon_1\varepsilon_2}\int_{\substack{|Q_1(\mathbf{y})|<\varepsilon_1\\|Q_2(\mathbf{y})|<\varepsilon_2}}U(\mathbf{y})\mathrm{d}\mathbf{y}=J_0(U)\ll H^{-n+2}.
$$ 
So it follows that
\begin{align*}
J_0(W)=J_0(\mathcal{P})+O(H^{-n+2}).
\end{align*}
This completes the proof of  Theorem~\ref{thm2}.\\


\end{document}